\documentclass[10pt,a4paper]{amsart}

\usepackage[T1]{fontenc}
\usepackage[utf8]{inputenc}
\usepackage{lmodern} 
\usepackage[english]{babel}
\usepackage{csquotes}
\usepackage[a4paper,margin=0.9in,marginparwidth=2cm]{geometry}
\usepackage{microtype}
\usepackage{xcolor}
\usepackage{graphicx}
\usepackage{amsmath,amssymb,amsthm,amsfonts,mathrsfs,amsopn}
\usepackage{mathtools}
\mathtoolsset{showonlyrefs} 
\usepackage{dsfont,esint,slashed}
\usepackage[all]{xy}
\usepackage{tikz-cd}
\usepackage{faktor}

\newif\ifdraftnotes
\draftnotesfalse
\newif\ifshowreferencekeys
\showreferencekeysfalse
\ifdraftnotes
  \usepackage[textsize=footnotesize]{todonotes}
\else
  \usepackage[disable]{todonotes}
\fi

\ifshowreferencekeys
  \usepackage[notref,notcite]{showkeys}
\fi

\usepackage[
  backend=biber,
  style=numeric-comp,
  sorting=nyt,
  sortcites=true,
  giveninits=true,
  maxbibnames=99,
  doi=true,
  url=true,
  eprint=true,
  isbn=false
]{biblatex}
\DeclareFieldFormat[article,inproceedings,unpublished,misc,online]{title}{#1}
\renewbibmacro{in:}{\ifentrytype{article}{}{\printtext{\bibstring{in}\intitlepunct}}}

\usepackage{xpatch}
\makeatletter
\xpatchbibmacro{cite:dump}
  {\ifnumgreater{\value{cbx@tempcnta}}{1}}
  {\ifnumgreater{\value{cbx@tempcnta}}{0}}
  {}
  {\PackageError{math-paper-template}{Citation compression patch failed}
    {Check the numeric-comp citation style in your TeX installation.}}
\makeatother

\AtEveryBibitem{
  \iffieldundef{doi}
    {\iffieldundef{eprint}
       {}
       {\clearfield{url}\clearfield{urlyear}\clearfield{urlmonth}\clearfield{urlday}}}
    {\clearfield{url}\clearfield{eprint}
     \clearfield{urlyear}\clearfield{urlmonth}\clearfield{urlday}}
}

\usepackage{xurl}
\usepackage{hyperref}
\definecolor{linkblue}{RGB}{25,55,105}
\hypersetup{
  colorlinks=true,
  linkcolor=linkblue,
  citecolor=linkblue,
  urlcolor=linkblue,
  pdfborder={0 0 0},
  pdftitle={Global solutions of the planar plasma equation with infinite p-mass},
  pdfauthor={Daniele Bartolucci, Aleks Jevnikar, Juncheng Wei, and Ruijun Wu}
}
\usepackage{bookmark}
\DeclareFieldFormat{doi}{%
  \mkbibacro{DOI}\addcolon\space
  \href{https://doi.org/#1}{\nolinkurl{#1}}}
\DeclareFieldFormat{eprint:arxiv}{%
  \newline\href{https://arxiv.org/abs/#1}{\nolinkurl{https://arxiv.org/abs/#1}}}
\DeclareFieldAlias{eprint:arXiv}{eprint:arxiv}
\DeclareFieldFormat{url}{\url{#1}}

\newcommand{\R}{\mathbb{R}}

\newcommand{\dd}{\mathop{}\!\mathrm{d}}

\newcommand{\sph}{\mathbb{S}}

\newcommand{\T}{\mathbb{T}}
\DeclareMathOperator{\ind}{ind}

\theoremstyle{plain}
\newtheorem{thm}{Theorem}[section]
\newtheorem{lemma}[thm]{Lemma}
\newtheorem{prop}[thm]{Proposition}

\newtheorem{theorem}{Theorem}
\newtheorem{mainprop}[theorem]{Proposition}

\theoremstyle{definition}

\theoremstyle{remark}
\newtheorem{rmk}[thm]{Remark}

\title[Global solutions of the planar plasma equation with infinite $p$-mass]{Global solutions of the planar plasma equation\\with infinite $p$-mass}
\date{}

\author[D. Bartolucci]{Daniele Bartolucci}
\address{Daniele Bartolucci, Department of Mathematics, University of Rome \textsl{"Tor Vergata"}, Via della ricerca scientifica n.1, 00133 Roma.}
\email{bartoluc@mat.uniroma2.it}

\author[A. Jevnikar]{Aleks Jevnikar}
\address{Aleks Jevnikar, Department of Mathematics, Computer Science and Physics, University of Udine, Via delle Scienze 206, 33100 Udine, Italy.}
\email{aleks.jevnikar@uniud.it}

\author[J. Wei]{Juncheng Wei}
\address{Juncheng Wei, Department of Mathematics, Chinese University of Hong Kong, Shatin, Hong Kong.}
\email{wei@math.cuhk.edu.hk}

\author[R. Wu]{Ruijun Wu}
\address{Ruijun Wu, School of mathematics and statistics, Beijing Institute of Technology, Zhongguancun South Street No. 5, 100081 Beijing, P.R.China.}
\email{ruijun.wu@bit.edu.cn}

\thanks{ D.B.
is partially supported by the MIUR Excellence Department Project MatMod@TOV
awarded to the Department of Mathematics, University of Rome ``Tor Vergata'' and by PRIN project 2026, ERC PE1\_11,
``{\em Variational Methods, Geometric Analysis and Free Boundary Problems}'', by INdAM-GNAMPA project ``{\em Analisi qualitativa di problemi differenziali non lineari}''
and by the E.P.G.P. Project sponsored by the University of Rome ``Tor Vergata''.\\
 \indent A.J. is partially supported by INdAM-GNAMPA project ``{\em Analisi qualitativa di problemi differenziali non lineari}''
 and PRIN Project 20227HX33Z ``{\em Pattern formation in nonlinear phenomena}''.\\
\indent J.W. is partially supported by GRF fund of RGC of Hong Kong entitled
``\emph{New frontiers in singularity formations of nonlinear partial differential equations}''.\\
\indent R.W. is supported by the Fundamental Research Funds for the Central Universities, the Beijing Institute of Technology Research Fund Program for Young Scholars, and the Beijing Natural Science Foundation No. 1262018.
}



\begin{document}

\begin{abstract}
We study entire solutions of~$-\Delta u=u_+^p$ in the plane for~$p>1$.
A solution which is nondecreasing in one direction is invariant in that direction, and for each nonconstant solution in this class the plasma region~$\{x:u(x)>0\}$ is a strip, the transverse profile has affine tails, and the Morse index is infinite.
The stable solutions are the nonpositive constants, and the nonconstant solutions with finite Morse index are precisely the translated and scaled radial profiles of finite~$p$-mass, each of which has Morse index one.
We also construct a continuous branch of solutions which are periodic in one variable, depend on both variables and have infinite Morse index.
\end{abstract}

\maketitle

\section{Introduction}

For~$p>1$, we consider the planar plasma equation
\begin{align}\label{eq:plasma}
 -\Delta u=u_+^p\quad\text{in }\R^2.
\end{align}
This equation arises in the analysis of singular limits of free boundary problems in plasma physics~\cite{BJWW}.
We call~$\{x:u(x)>0\}$ the plasma region of~$u$, and the~$p$-mass of~$u$ is
\begin{align}
 M_p(u)=\int_{\R^2}u_+^p\dd x .
\end{align}
Every nonconstant solution of finite~$p$-mass has a disk as its plasma region and is determined, up to translation and dilation, by the positive radial Dirichlet solution on the unit disk~\cite[Proposition~2.2]{BJWW}.
There are also solutions of infinite~$p$-mass depending on one linear coordinate, for which the plasma region is a strip~\cite[Section~3.2]{BJWW}.
We ask when a global solution must have one of these forms and whether there are solutions of infinite~$p$-mass which depend on both coordinates.

We begin with solutions which are monotone in some direction.
The planar case of the De Giorgi conjecture~\cite{DeGiorgi1979,GhoussoubGui1998} concerns bounded entire solutions of the Allen--Cahn equation
\begin{align}
 -\Delta u=u-u^3 \quad\text{in }\R^2
\end{align}
for which~$\partial_eu>0$ in some direction~$e$.
The level sets of such a solution are parallel lines, and for some unit vector~$\nu$ and some~$b\in\R$,
\begin{align}
 u(x)=H(\nu\cdot x-b),\qquad H(t)=\tanh(t/\sqrt{2}).
\end{align}
Differentiating the equation for the profile~$H$, we find that its derivative defines a positive solution of the planar linearized equation,
\begin{align}
 h(x)=H'(\nu\cdot x-b)>0,\qquad -\Delta h=(1-3u^2)h .
\end{align}
For each~$\xi\in C_c^1(\R^2)$, testing this equation with~$\xi^2/h$ and completing the square proves stability, since
\begin{align}
 \int_{\R^2}\bigl(|\nabla\xi|^2-(1-3u^2)\xi^2\bigr)\dd x
 =\int_{\R^2}h^2\left|\nabla\!\left(\frac{\xi}{h}\right)\right|^2\dd x\ge0.
\end{align}
For the three-dimensional symmetry results, see~\cite{AmbrosioCabre2000,AlbertiAmbrosioCabre2001}.
Related symmetry and Liouville theorems for entire solutions are proved in~\cite{Farina2003,Farina2005,FarinaValdinoci}, while~\cite{BCN1997,FarinaSciunzi2016,FarinaSciunzi2017,BeuvinFarina2026} treat unbounded domains, including half-spaces and strips.

For the plasma equation we impose the weak inequality~$\partial_e u\ge0$ and allow~$u$ to be unbounded.
Here the direction in the hypothesis is itself a direction of invariance.

\begin{theorem}\label{thm:rigidity}
Let~$g\in C(\R)$ be monotone and suppose that~$u\in C^2(\R^2)$ satisfies
\begin{align}
 -\Delta u=g(u),\qquad \partial_e u\ge0\quad\text{in }\R^2
\end{align}
for some unit vector~$e$.
Fix a unit vector~$\nu\perp e$ and write~$x=se+t\nu$.
\begin{enumerate}
\renewcommand{\labelenumi}{(\roman{enumi})}
\item If~$g\equiv c$, then for some~$d\ge0$ and~$a,b\in\R$,
\begin{align}
 u(se+t\nu)=ds-\tfrac{c}{2}t^2+at+b.
\end{align}
\item If~$g$ is nondecreasing and nonconstant, then~$\partial_e u\equiv0$.
If moreover~$g\ge0$ and~$u$ is nonconstant, then for some~$b\in\R$ and~$a$ with~$g(a)>0$,
\begin{align}
 u(x)&=U_a(\nu\cdot x-b),\\
 -U_a''&=g(U_a),\qquad U_a(0)=a,\quad U_a'(0)=0.
\end{align}
For each~$a$ with~$g(a)>0$, the stated initial conditions determine a unique global profile~$U_a$.
It is even, strictly decreasing on~$(0,\infty)$ and tends to~$-\infty$ at both ends.
\item If~$g$ is nonincreasing and nonconstant, then either~$\partial_e u\equiv0$ or~$\partial_e u$ is unbounded above.
\end{enumerate}
\end{theorem}

For~\eqref{eq:plasma} we write~$W_a=U_a$, where~$a>0$.
The positive part of~$W_a$ is supported on an interval~$[-R_a,R_a]$, beyond which~$W_a$ is affine on each half-line.
In particular every nonconstant solution in Theorem~\ref{thm:rigidity} for the plasma equation has infinite~$p$-mass.

The profiles~$W_a$ are invariant in the direction~$e$ from the hypothesis, so~$\partial_e u$ vanishes identically.
Their stability must be decided from the second variation.
For~$g\in C^1(\R)$ the quadratic form at an entire solution of~$-\Delta u=g(u)$ is
\begin{align}\label{eq:form}
 Q_u(\xi)=\int_{\R^2}\bigl(|\nabla\xi|^2-g'(u)\xi^2\bigr)\dd x,\qquad \xi\in C_c^1(\R^2).
\end{align}
The Morse index of~$u$, denoted by~$\ind(u)$, is the supremum of the dimensions of subspaces of~$C_c^1(\R^2)$ on which~$Q_u$ is negative definite.
Thus~$u$ is stable if and only if~$\ind(u)=0$.
For interior regularity estimates for stable solutions in dimensions up to nine, see~\cite[Theorem~1.2]{CFRS}.
In the plane the sign~$g'\ge0$ allows us to test stability with logarithmic cutoffs of arbitrarily small Dirichlet energy, from which the classification of stable solutions follows without a growth assumption at infinity.

Let~$\phi$ denote the unique positive radial solution of
\begin{align}
\begin{cases}
     -\Delta\phi=\phi^p &\text{in }B_1,\\
 \phi=0 &\text{on }\partial B_1.
\end{cases}
\end{align}
Here and below~$B_R(x_0)$ denotes the open ball of radius~$R$ centered at~$x_0$, and we write~$B_R=B_R(0)$.
For~$R>0$ and~$z\in\R^2$, let
\begin{align}\label{eq:diskprofile}
 U_{R,z}(x)=
 \begin{cases}
  R^{-2/(p-1)}\phi(r/R),&0\le r\le R,\\
  R^{-2/(p-1)}\phi'(1)\log(r/R),&r>R,
 \end{cases}
 \qquad r=|x-z|.
\end{align}

The stable solutions of~\eqref{eq:plasma} are precisely the nonpositive constants, and every solution with finite Morse index has compactly supported positive part~\cite[Theorems~1.3 and~1.5]{RW}.
In particular such a solution has finite~$p$-mass, so by~\cite[Proposition~2.2]{BJWW} every nonconstant function in this class is one of the profiles in~\eqref{eq:diskprofile}.
We compute its index with respect to all compactly supported variations in the plane.

\begin{mainprop}\label{thm:index}
For every~$R>0$ and~$z\in\R^2$, the radial solution~$U_{R,z}$ has Morse index one on~$\R^2$.
\end{mainprop}

It also follows from~\cite[Theorems~1.3 and~1.5]{RW} that every nonconstant one-dimensional solution of~\eqref{eq:plasma} has infinite Morse index in the plane.
For comparison, four-ended planar Allen--Cahn solutions of index one were constructed in~\cite{GLW}.

We turn to semi-periodic solutions, motivated by the pioneering work of Dancer~\cite{Dancer2001}, who constructed solutions bifurcating from a decaying profile and periodic in one additional variable.
For the planar plasma equation we use the strip profile~$W_1$, that is,~$W_a$ with~$a=1$, so that~$W_1(0)=1$.
The solutions on the branch are periodic in~$y$ and their plasma regions lie in a fixed strip.
By Proposition~\ref{prop:profile}, the transverse operator on even functions in~$L^2(\R)$ has a lowest eigenvalue~$\lambda_0<0$ with a positive normalized eigenfunction~$\psi_0$, so that
\begin{align}
 -\psi_0''-p(W_1)_+^{p-1}\psi_0=\lambda_0\psi_0,\qquad \omega_0=\sqrt{-\lambda_0}.
\end{align}

The space~$\mathscr{X}$ in the following theorem consists of perturbations~$v(x,t)$ which are even in both variables and~$2\pi$-periodic in~$t$.
Writing
\begin{align}
 v(x,t)=\bar v(x)+w(x,t),\qquad
 \bar v(x)=\frac{1}{2\pi}\int_0^{2\pi}v(x,t)\dd t ,
\end{align}
we require the mean~$\bar v$ to vanish at the origin and allow it to be asymptotically affine on each half-line, while the component~$w$ of zero mean decays exponentially as~$|x|\to\infty$.
The precise definition and norm are given at the beginning of Section~\ref{sec:unified}.
\begin{theorem}\label{thm:main}
For each~$p>1$ there is a continuous branch~$(\omega(s),v_s)$ defined for~$|s|<\epsilon$ for some~$\epsilon>0$.
The branch passes through~$(\omega_0,0)$ at~$s=0$, and the function
\begin{align}
 u_s(x,y)=W_1(x)+v_s(x,\omega(s)y)
\end{align}
is a classical solution of~\eqref{eq:plasma}, with periodic perturbation satisfying
\begin{align}\label{eq:planarbranch}
 v_s(x,t)=s\psi_0(x)\cos t+O_{\mathscr{X}}\bigl(|s|^{\min\{p,2\}}\bigr).
\end{align}
For~$s\ne0$ sufficiently small,~$u_s$ depends on both linear coordinates, its positive part is uniformly bounded and supported in a fixed strip, and both its~$p$-mass and its Morse index on~$\R^2$ are infinite.
\end{theorem}

\begin{rmk}\label{rem:regularity}
If~$p>2$, the branch~$s\mapsto(\omega(s),v_s)$ is of class~$C^1$ in~$\R\times\mathscr{X}$.
\end{rmk}

For the plasma profile, the linearized potential vanishes outside the positive strip, where the zero Fourier mode has affine tails.
We allow those tails and fix the mean at the origin, while the nonzero Fourier modes are controlled in an exponentially weighted space.
This leads to a simple eigenvalue bifurcation in which the divided equation is continuous in the amplitude for every~$p>1$.

We prove Theorem~\ref{thm:rigidity} in Section~\ref{sec:monotone}, first establishing invariance in the direction of monotonicity and then determining the transverse profile.
The proof of Proposition~\ref{thm:index} is in Subsection~\ref{subsec:finite-index}, where we pass from exterior stability to finite~$p$-mass and compute the Morse index of the radial solutions.
Theorem~\ref{thm:main} follows at the end of Subsection~\ref{subsec:positivity-index} by taking~$k=l=1$ in Theorem~\ref{thm:semi}, after constructing the semi-periodic branch in Subsection~\ref{subsec:bifurcation}.

\section{Solutions monotone in a direction}\label{sec:monotone}

\begin{proof}[Proof of Theorem~\ref{thm:rigidity}]
For each~$h>0$ the increment~$v_h(x)=u(x+he)-u(x)$ is nonnegative by directional monotonicity, and
\begin{align}
 -\Delta v_h(x)=g(u(x+he))-g(u(x)).
\end{align}
We apply the planar Liouville theorem in the following cases.
\begin{itemize}
\item If~$g$ is constant or nondecreasing, then~$v_h$ is superharmonic and bounded below, hence constant.
\item If~$g$ is nonincreasing and~$\partial_e u$ is bounded above by~$M$, then~$v_h$ is subharmonic and bounded above by~$Mh$, hence constant.
\end{itemize}
In either case~$v_h(x)=v_h(0)$.
Dividing by~$h$ and letting~$h\downarrow0$, we obtain
\begin{align}
 \partial_e u(x)=\partial_e u(0)=:d\ge0.
\end{align}

Thus, with~$U(t)=u(t\nu)$,
\begin{align}
 u(se+t\nu)&=ds+U(t),\\
 -U''(t)&=g(ds+U(t))\qquad(s,t\in\R).
\end{align}
If~$g\equiv c$, the first assertion follows by integrating the second identity twice.
If~$g$ is nonconstant, then~$d=0$, since for~$d>0$ the argument on the right ranges over~$\R$ with~$t$ fixed.
This proves the directional assertions, including the alternative for nonincreasing~$g$.

It remains to classify the nonconstant profiles when~$g$ is nonnegative, nondecreasing and nonconstant.
Since~$U''=-g(U)\le0$, the derivative~$U'$ is nonincreasing.
There is a point~$t_0$ with~$g(U(t_0))>0$, since otherwise~$U$ would be a nonconstant affine function with range~$\R$ and hence~$g$ would vanish on~$\R$.
We claim that~$U$ attains its maximum.
If~$U'(t_0)>0$, then as long as~$U'$ stays positive to the right of~$t_0$, we have~$U(t)\ge U(t_0)$ and
\begin{align}
 U'(t)=U'(t_0)-\int_{t_0}^t g(U(s))\dd s
       \le U'(t_0)-(t-t_0)g(U(t_0)).
\end{align}
The right-hand side reaches zero in finite time, so~$U'$ has a zero.
If~$U'(t_0)<0$, the same argument applies to~$t\mapsto U(-t)$, while for~$U'(t_0)=0$ we already have a critical point.
Each zero of~$U'$ is a global maximum by concavity.
Translate one such point to~$0$ and put~$a=U(0)$.
Then~$g(a)\ge g(U(t_0))>0$, so~$U''(0)<0$.
Since~$U'$ is nonincreasing, it is strictly positive on~$(-\infty,0)$ and strictly negative on~$(0,\infty)$.

To determine both branches, set
\begin{align}
 H_a(s)=\int_s^a g(r)\dd r,\qquad s\le a.
\end{align}
Multiplying the equation by~$U'$ and integrating from~$0$ to~$t$, we obtain
\begin{align}
 \frac12 U'(t)^2
   =-\int_0^t g(U(\tau))U'(\tau)\dd\tau
   =\int_{U(t)}^a g(r)\dd r
   =H_a(U(t)).
\end{align}
Since~$g(a)>0$, we have~$H_a(s)>0$ for~$s<a$.
Taking the sign of~$U'$ on each half-line into account, we have
\begin{align}
 U'(t)=
 \begin{cases}
  -\sqrt{2H_a(U(t))},&t>0,\\
   \sqrt{2H_a(U(t))},&t<0.
 \end{cases}
\end{align}
Separating variables on the two half-lines, we arrive at the same formula,
\begin{align}\label{eq:quadrature}
 |t|=\int_{U(t)}^a\frac{\dd s}{\sqrt{2H_a(s)}} .
\end{align}
For~$t>0$, applying~\eqref{eq:quadrature} at~$t$ and~$-t$, we find
\begin{align}
 \int_{U(t)}^a\frac{\dd s}{\sqrt{2H_a(s)}}
 =t=\int_{U(-t)}^a\frac{\dd s}{\sqrt{2H_a(s)}}.
\end{align}
Since the integrand is strictly positive, equality of the two integrals forces their lower endpoints to coincide, so~$U$ is even.

It remains to show that each~$a$ with~$g(a)>0$ is the maximum of a unique profile centered at the origin.
As~$s\uparrow a$, continuity implies
\begin{align}
 H_a(s)=g(a)(a-s)+o(a-s),
\end{align}
 so the singularity at the upper endpoint of~\eqref{eq:quadrature} is integrable.
By monotonicity we have
\begin{align}
 H_a(s)\le g(a)(a-s)\qquad(s<a),
\end{align}
from which it follows that
\begin{align}
 \int_s^a\frac{\dd r}{\sqrt{2H_a(r)}}\ge
 \int_s^a\frac{\dd r}{\sqrt{2g(a)(a-r)}}=\sqrt{\frac{2(a-s)}{g(a)}} .
\end{align}
Thus the integral in~\eqref{eq:quadrature} ranges from~$0$ to~$+\infty$ as its lower endpoint decreases from~$a$ to~$-\infty$.
Its inverse defines the unique decreasing branch~$U_a(t)$ for~$t\ge0$, with~$U_a(t)\to-\infty$ as~$t\to\infty$.
For~$t>0$, we differentiate the inverse relation and use~$H_a'=-g$ to obtain
\begin{align}
 U_a'&=-\sqrt{2H_a(U_a)},\\
 U_a''&=-\frac{H_a'(U_a)U_a'}{\sqrt{2H_a(U_a)}}=-g(U_a).
\end{align}
As~$t\downarrow0$, these derivatives tend to~$0$ and~$-g(a)$, respectively, so the even extension is~$C^2$ across the origin.
Restoring the position~$b$ of the maximum, we have
\begin{align}
 u(x)=U_a(\nu\cdot x-b).
\end{align}

\end{proof}

For nonincreasing~$g$ as in Theorem~\ref{thm:rigidity}~(iii), the following oscillation estimate holds for every entire solution, even without directional monotonicity.
For every~$x_0\in\R^2$ and~$R>0$, we have that
\begin{align}
 \operatorname{osc}_{B_R(x_0)}u
 =\sup_{B_R(x_0)}u-\inf_{B_R(x_0)}u
 \ge \frac{R^2}{4}\,|g(u(x_0))|.
\end{align}
Indeed, if~$g(u(x_0))>0$, fix~$0<a<g(u(x_0))$ and consider
\begin{align}
 w(x)=u(x)+\tfrac{a}{4}|x-x_0|^2 .
\end{align}
At an interior minimum~$y$ of~$w$ on the closed disk we would have~$w(y)\le w(x_0)$ and thus~$u(y)\le u(x_0)$.
Since~$g$ is nonincreasing, this would imply
\begin{align}
 0\le\Delta w(y)=-g(u(y))+a\le -g(u(x_0))+a<0.
\end{align}
Hence the minimum is attained on the boundary, where the quadratic term equals~$aR^2/4$, and
\begin{align}
 \min_{\partial B_R(x_0)}u+\frac{aR^2}{4}
 =\min_{\overline{B_R(x_0)}}w
 \le w(x_0)=u(x_0).
\end{align}
By continuity the oscillation on~$B_R(x_0)$ is at least the difference between~$u(x_0)$ and the boundary minimum, so we may let~$a\uparrow g(u(x_0))$ to obtain the estimate.
If~$g(u(x_0))<0$, apply the same argument to~$-u$ with the nonincreasing nonlinearity~$s\mapsto-g(-s)$.
For a nonconstant~$u$ there is some~$x_0$ with~$g(u(x_0))\ne0$, since otherwise~$u$ would be a nonconstant entire harmonic function and hence take every real value, forcing~$g\equiv0$.
Thus every nonconstant solution has at least quadratic oscillation on large disks.

To see that a solution in case~(iii) may depend on both coordinates, take~$g(s)=(-s)_+$ and
\begin{align}
 u(se+t\nu)=-e^{-s}-e^{-t}.
\end{align}
Indeed,
\begin{align}
 -\Delta u=-u=g(u),\qquad \partial_e u=e^{-s},
\end{align}
so the directional derivative is unbounded above.
A direction of invariance would require
\begin{align}
 \alpha e^{-s}+\beta e^{-t}=0
\end{align}
 for all~$s,t$, which forces~$\alpha=\beta=0$.

For the plasma equation, the profiles~$W_a$ in Theorem~\ref{thm:rigidity}~(ii) are precisely those constructed in~\cite[Section~3.2]{BJWW}.
Their plasma regions are strips, and each profile has infinite~$p$-mass in the plane.

\section{Stable solutions and solutions with finite Morse index}\label{sec:index}

\subsection{Global stability}
We first determine what global stability implies for the more general nonlinearity in Theorem~\ref{thm:rigidity}.
The next lemma is a consequence of the spectral result in~\cite[Section~1]{Shargorodsky2014}, applied to a bounded, compactly supported truncation of the potential.
For completeness, we include a direct proof.

\begin{lemma}\label{lemma:capacity}
Let~$V\ge0$ be locally integrable on~$\R^2$ and suppose that, for every~$\xi\in C_c^1(\R^2)$,
\begin{align}
 \int_{\R^2} V\xi^2\dd x\le\int_{\R^2}|\nabla\xi|^2\dd x.
\end{align}
Then~$V=0$ almost everywhere.
\end{lemma}
\begin{proof}
Fix~$0<\rho<R$ and take the radial cutoff~$\eta$ which equals one on~$B_\rho$, vanishes outside~$B_R$, and is linear in~$\log r$ on the intervening annulus.
Approximating~$\eta$ by smooth cutoffs uniformly and in~$H^1$, we obtain
\begin{align}
 0\le\int_{B_\rho}V\dd x
 &\le\int_{\R^2}|\nabla\eta|^2\dd x
 =\frac{2\pi}{\log^2(R/\rho)}\int_\rho^R\frac{\dd r}{r}
 =\frac{2\pi}{\log(R/\rho)}.
\end{align}
Letting~$R\to\infty$ proves the claim for every~$\rho>0$.
\end{proof}

\begin{prop}\label{prop:plateau}
Suppose that~$g\in C^1(\R)$ is nonnegative, nondecreasing and nonconstant.
Let~$u$ be a stable solution of~$-\Delta u=g(u)$ in~$\R^2$.
Then either~$u$ is constant, or
\begin{align}\label{eq:plateau}
 u(x)=a-\frac12(x-z)^T A(x-z)
\end{align}
for some~$a\in\R$,~$z\in\R^2$ and nonzero symmetric positive semidefinite matrix~$A$ such that
\begin{align}
 g(s)=\operatorname{tr} A\qquad\text{for every }s\le a.
\end{align}
In particular, suppose that there are no~$a\in\R$ and~$c>0$ such that~$g(s)=c$ for every~$s<a$.
Then every stable solution is constant, and every nonconstant one-dimensional solution has infinite Morse index in the plane.
\end{prop}
\begin{proof}
Since~$g'(u)\ge0$, stability and Lemma~\ref{lemma:capacity} imply that~$g'(u)=0$ everywhere, so~$g$ is constant on the interval~$u(\R^2)$, say~$g(u)=c$.
Suppose~$u$ is nonconstant.
It is unbounded below by planar superharmonicity, as in Section~\ref{sec:monotone}.
Its supremum~$a$ must be finite, since otherwise its range would be~$\R$ and~$g$ would be constant on~$\R$.
If~$c=0$, the function is harmonic and bounded above, hence constant.
Thus~$c>0$, and the function~$h=u+c|x|^2/4$ is harmonic and bounded above by~$a+c|x|^2/4$.
For~$y\in\R^2$ and~$R>0$, set
\begin{align}
 M_R=a+\frac{c}{4}(|y|+2R)^2+1.
\end{align}
The function~$M_R-h$ is positive and harmonic on~$B_{2R}(y)$.
By the interior derivative estimate on~$B_R(y)$ and Harnack's inequality on~$B_{2R}(y)$, we have
\begin{align}
 |D^3h(y)|&\le CR^{-3}\sup_{B_R(y)}(M_R-h)\\
 &\le CR^{-3}(M_R-h(y))\\
 &\le CR^{-3}\bigl(a+1+c(|y|+2R)^2/4-h(y)\bigr).
\end{align}
For fixed~$y$ the last expression tends to zero as~$R\to\infty$, so~$D^3h=0$ throughout the plane.
Thus~$u$ is a quadratic polynomial with Hessian~$-A$.
Boundedness above forces~$A\ge0$ and the linear coefficient to vanish on~$\ker A$, hence this coefficient is~$Az$ for some~$z$.
Completing the square, we obtain~\eqref{eq:plateau}, whose range is~$(-\infty,a]$.
Since~$-\Delta u=\operatorname{tr} A=c>0$, the matrix~$A$ is nonzero and~$g(s)=\operatorname{tr} A$ for every~$s\le a$.
Conversely,~$g'(u)=0$ for these quadratics and for the constants at zeros of the nonnegative function~$g$, so they are stable.

Under the additional hypothesis on~$g$, the quadratic case is excluded.
For any nonconstant one-dimensional solution, there is then a compactly supported~$\xi$ with~$Q_u(\xi)<0$.
Translations of~$\xi$ in the invariant direction have the same quadratic form, and arbitrarily many such translates have disjoint supports.
Their span is negative definite, which proves the last assertion.
\end{proof}

For example,~$g(s)=1+s_+^2$ admits the stable solutions~$a-x_1^2/2$ and~$a-|x|^2/4$ for every~$a\le0$.
The first depends on one coordinate, while the second has a Hessian of rank two.
For~$g(s)=s_+^p$ the only stable solutions are~$u\equiv c\le0$.

\subsection{Finite Morse index for the plasma equation}\label{subsec:finite-index}
We now specialize to the power nonlinearity.
The first step connects stability outside a compact set with finite~$p$-mass and uses only tests in the region where stability holds.

\begin{lemma}\label{lemma:annuli}
If a solution~$u$ of~\eqref{eq:plasma} is stable outside~$B_{R_0}(x_0)$, then for every~$R>2R_0$ we have
\begin{align}\label{eq:annuli}
 \begin{aligned}
 \int_{B_{2R}(x_0)\setminus B_R(x_0)}u_+^{p+1}\dd x &\le C_p R^{-4/(p-1)},\\
 \int_{\R^2\setminus B_R(x_0)}u_+^p\dd x &\le C_p R^{-2/(p-1)}.
 \end{aligned}
\end{align}
\end{lemma}
\begin{proof}
Put~$w=u_+$ and take~$\eta\in C_c^1(\R^2\setminus\overline{B_{R_0}(x_0)})$.
Since~$\nabla w=\mathbf{1}_{\{u>0\}}\nabla u$ almost everywhere, we may test the equation with~$w\eta^2$ to obtain
\begin{align}
 \int_{\R^2} w^{p+1}\eta^2\dd x
 &=\int_{\R^2}\eta^2|\nabla w|^2\dd x
 +2\int_{\R^2}w\eta\,\nabla w\cdot\nabla\eta\dd x.
\end{align}
Testing the stability inequality with~$w\eta$ and expanding the gradient, we find
\begin{align}
 p\int_{\R^2} w^{p+1}\eta^2\dd x
 \le\int_{\R^2}|\nabla(w\eta)|^2\dd x
 =\int_{\R^2}\eta^2|\nabla w|^2\dd x
 +2\int_{\R^2}w\eta\,\nabla w\cdot\nabla\eta\dd x
 +\int_{\R^2}w^2|\nabla\eta|^2\dd x.
\end{align}
Both tests follow by~$H^1$ approximation on their compact supports.
Subtracting the first identity from the last inequality, we obtain
\begin{align}
 (p-1)\int_{\R^2} w^{p+1}\eta^2\dd x
 \le\int_{\R^2}w^2|\nabla\eta|^2\dd x.
\end{align}
Write~$A_R=B_{2R}(x_0)\setminus B_R(x_0)$ and take~$0\le\chi\le1$ supported in~$B_{4R}(x_0)\setminus\overline{B_{R/2}(x_0)}$, equal to one on~$A_R$, with~$|\nabla\chi|\le C/R$.
Since~$R>2R_0$ this cutoff is supported in the stable region.
Choose the least integer~$m\ge(p+1)/(p-1)$, take~$\eta=\chi^m$ and set
\begin{align}
 J=\int_{\R^2}w^{p+1}\chi^{2m}\dd x.
\end{align}
Since~$2m-2\ge4m/(p+1)$, we can apply H\"older's inequality in the form
\begin{align}
 J&\le C_p R^{-2}\int_{B_{4R}(x_0)} w^2\chi^{2m-2}\dd x\\
 &\le C_p R^{-2}\int_{B_{4R}(x_0)}(w^{p+1}\chi^{2m})^{2/(p+1)}\dd x\\
 &\le C_p R^{-2}|B_{4R}(x_0)|^{(p-1)/(p+1)}J^{2/(p+1)}
 \le C'_p R^{-4/(p+1)}J^{2/(p+1)}.
\end{align}
For~$J>0$ the first estimate follows on dividing by~$J^{2/(p+1)}$, and for~$J=0$ it is immediate.
Applying H\"older's inequality on~$A_R$ and then summing over dyadic annuli, we find
\begin{align}
 \int_{A_R}w^p\dd x
 &\le \left(\int_{A_R}w^{p+1}\dd x\right)^{p/(p+1)}|A_R|^{1/(p+1)}
 \le C_pR^{-2/(p-1)},\\
 \int_{\R^2\setminus B_R(x_0)}w^p\dd x
 &=\sum_{j=0}^{\infty}\int_{A_{2^jR}}w^p\dd x
 \le C_pR^{-2/(p-1)}\sum_{j=0}^{\infty}2^{-2j/(p-1)}.
\end{align}
The series converges since~$p>1$.
\end{proof}

\begin{proof}[Proof of Proposition~\ref{thm:index}]
Translation and a change of scale preserve the index, so we take~$U=U_{1,0}$ and set~$V=pU_+^{p-1}$.
By~\eqref{eq:diskprofile}, the function~$U$ agrees with~$\phi$ on~$\overline{B_1}$ and has a logarithmic tail outside the unit disk.
Testing the equation for~$U$ with~$U_+$, we find
\begin{align}
 Q_U(U_+)=\int_{B_1}\bigl(|\nabla\phi|^2-p\phi^{p+1}\bigr)\dd x
 =(1-p)\int_{B_1}\phi^{p+1}\dd x<0.
\end{align}
Smooth approximation of~$U_+$ in~$H^1(\R^2)$ produces a compactly supported smooth test with negative quadratic form.

For the upper bound, take a real-valued test function~$\xi\in C_c^\infty(\R^2)$ and write its angular Fourier expansion in polar coordinates,
\begin{align}
 \xi(r,\theta)
 &=\xi_0(r)+\sum_{j=1}^{\infty}\bigl(\xi_j^c(r)\cos(j\theta)+\xi_j^s(r)\sin(j\theta)\bigr),\\
 \xi_0(r)&=\frac{1}{2\pi}\int_0^{2\pi}\xi(r,\theta)\dd\theta.
\end{align}
For a radial coefficient~$\eta$ of frequency~$j$, write
\begin{align}
 q_j[\eta]=\int_0^\infty\left(|\eta'|^2+
 \left(\frac{j^2}{r^2}-V(r)\right)\eta^2\right)r\dd r.
\end{align}
Since~$V$ is radial, orthogonality eliminates the terms involving distinct angular modes, while differentiating in~$\theta$ supplies the factor~$j^2$.
By Parseval's identity,
\begin{align}
 Q_U(\xi)
 &=\int_0^\infty\int_0^{2\pi}\left(|\partial_r\xi|^2+
 \frac{1}{r^2}|\partial_\theta\xi|^2-V(r)\xi^2\right)r\dd\theta\dd r\\
 &=2\pi q_0[\xi_0]+\pi\sum_{j=1}^{\infty}\bigl(q_j[\xi_j^c]+q_j[\xi_j^s]\bigr).
\end{align}

For the nonzero modes we use the radial factor~$h=-U'$ of the translation derivative~$-\partial_{x_1}U=h(r)\cos\theta$.
The radial equation for~$U$ and its integral from the origin read
\begin{align}
 h'(r)+\frac{h(r)}r&=U_+(r)^p, &
 rh(r)&=\int_0^r tU_+(t)^p\dd t,
\end{align}
so~$h(r)>0$ for every~$r>0$.
Since~$s\mapsto s_+^p$ is~$C^1$, differentiation of the first identity is valid across~$r=1$ and implies
\begin{align}
 h''+\frac1r h'-\frac1{r^2}h
 &=pU_+^{p-1}U'=-Vh,\\
 (rh')'&=r\left(\frac1{r^2}-V\right)h .
\end{align}
To determine the boundary behavior, use~$h(0)=0$ and substitute~$t=rs$ in the integral for~$rh(r)$,
\begin{align}
 h'(0)=-\phi''(0)=\lim_{r\downarrow0}\frac{h(r)}r
 =\lim_{r\downarrow0}\int_0^1sU_+(rs)^p\dd s=\frac{\phi(0)^p}{2}.
\end{align}
Since~$h'(0)>0$, it follows that~$rh'(r)/h(r)\to1$ as~$r\downarrow0$.

Take~$j\ge1$ and let~$b_j(\theta)$ be either~$\cos(j\theta)$ or~$\sin(j\theta)$, with corresponding coefficient~$\eta$.
As~$b_j$ has mean zero, we have
\begin{align}
 \eta(r)=\frac1\pi\int_0^{2\pi}\bigl(\xi(r,\theta)-\xi(0)\bigr)b_j(\theta)\dd\theta=O(r),
\end{align}
where smoothness of~$\xi$ bounds the difference in parentheses uniformly in~$\theta$.
Expanding the square below and using the equation for~$h$, we obtain
\begin{align}
 r\left(|\eta'|^2+\left(\frac1{r^2}-V(r)\right)\eta^2\right)
 =rh^2\left|\left(\frac{\eta}{h}\right)'\right|^2
 +\left(r\frac{h'}{h}\eta^2\right)'.
\end{align}
Since~$\eta(r)=O(r)$ and~$rh'(r)/h(r)\to1$, the boundary term vanishes at the origin, and it vanishes at infinity by compact support.
Integration and addition of the remaining angular term now show that
\begin{align}
 q_j[\eta]
 =\int_0^\infty rh^2\left|\left(\frac{\eta}{h}\right)'\right|^2\dd r
 +(j^2-1)\int_0^\infty\frac{\eta^2}{r}\dd r\ge0 .
\end{align}

For the circular mean~$\xi_0$, we obtain a radial solution of the linearized equation by varying the radius in the family~$U_{R,0}$.
By~\eqref{eq:diskprofile},
\begin{align}
 U_{R,0}(r)&=R^{-2/(p-1)}U(r/R),\\
 Z(r)&=-\left.\frac{\partial}{\partial R}U_{R,0}(r)\right|_{R=1}
 =\frac{2}{p-1}U(r)+rU'(r).
\end{align}
Each profile solves the plasma equation, so differentiation in~$R$ implies
\begin{align}
 (rZ')'=-rVZ.
\end{align}
Using the radial equation inside the unit disk and the logarithmic formula outside, we obtain
\begin{align}
 Z'&=\left(1+\frac{2}{p-1}\right)U'+rU''
     =\frac{2}{p-1}U'-rU^p<0 &&(0<r<1),\\
 Z&=\phi'(1)\left(1+\frac{2}{p-1}\log r\right)<0 &&(r\ge1).
\end{align}
Since~$Z(0)>0$ and~$Z(1)=\phi'(1)<0$, the function~$Z$ has exactly one zero~$r_*\in(0,1)$, which is simple.
Suppose~$\xi_0(r_*)=0$.
On each of the intervals~$(0,r_*)$ and~$(r_*,\infty)$, the equation for~$Z$ implies
\begin{align}
 r\bigl(|\xi_0'|^2-V\xi_0^2\bigr)
 =rZ^2\left|\left(\frac{\xi_0}{Z}\right)'\right|^2
 +\left(r\frac{Z'}{Z}\xi_0^2\right)'.
\end{align}
The boundary term tends to zero on either side of~$r_*$ since~$Z$ has a simple zero there and~$\xi_0(r)=O(|r-r_*|)$.
It also vanishes at the origin, since~$Z(0)>0$ and~$Z'(0)=0$, and at infinity by compact support.
Integrating separately on the two intervals, we obtain
\begin{align}
 q_0[\xi_0]
 =\int_0^{r_*}rZ^2\left|\left(\frac{\xi_0}{Z}\right)'\right|^2\dd r
 +\int_{r_*}^{\infty}rZ^2\left|\left(\frac{\xi_0}{Z}\right)'\right|^2\dd r\ge0 .
\end{align}
Combining the angular modes, we have proved that~$Q_U(\xi)\ge0$ whenever~$\xi_0(r_*)=0$.
On any two-dimensional space of smooth tests, the linear functional~$\xi\mapsto\xi_0(r_*)$ has a nonzero function in its kernel.
Hence the quadratic form cannot be negative definite on a two-dimensional space of smooth tests.
The same upper bound on the index holds on~$C_c^1(\R^2)$ by smooth approximation of a basis of any finite-dimensional negative definite subspace.
Together with the negative direction obtained from~$U_+$, this proves Proposition~\ref{thm:index}.
\end{proof}

By~\cite[Theorems~1.3 and~1.5]{RW}, stability forces a nonpositive constant and finite index forces compact support of~$u_+$.
Lemma~\ref{lemma:annuli} also shows directly that exterior stability implies finite~$p$-mass.
By~\cite[Proposition~2.2]{BJWW} a solution of finite~$p$-mass is either a nonpositive constant or one of the profiles~\eqref{eq:diskprofile}.
The former have index zero, while the latter have index one by Proposition~\ref{thm:index}.

\subsection{Other cutoff nonlinearities}
The existence of a radial solution of finite positive index extends to a general cutoff, even when the nonlinearity has no scaling symmetry.

\begin{prop}\label{prop:generalcore}
Suppose that~$g\in C^1(\R)$ is nonnegative, nondecreasing and nonconstant, and that~$g$ has at least one zero.
Let~$b$ be the largest zero of~$g$.
For each~$a>b$ there is a radial entire solution of~$-\Delta u=g(u)$ with~$u(0)=a$, of finite positive Morse index, which is strictly greater than~$b$ on a disk and has a logarithmic tail outside that disk.
\end{prop}
\begin{proof}
We have~$g=0$ on~$(-\infty,b]$ and~$g>0$ on~$(b,\infty)$.
The regular radial initial-value problem has~$u'(0)=0$ and
\begin{align}
 -ru'(r)=\int_0^r t g(u(t))\dd t.
\end{align}
Since~$u\le a$, the right-hand side lies between zero and~$g(a)r^2/2$, so the solution continues globally.
It reaches~$b$ at a finite radius~$R$, since the integral is bounded below by a positive constant for large~$r$ and
\begin{align}
 \int_1^\infty\frac{\dd r}{r}=\infty.
\end{align}
For~$r\ge R$ it equals
\begin{align}
 u(r)=b-m\log(r/R),\qquad m=-Ru'(R)>0.
\end{align}
The potential~$g'(u)$ is bounded, compactly supported and nonzero because~$g(a)>g(b)$.
By Lemma~\ref{lemma:capacity}, the quadratic form has a negative direction.
For an upper bound on the index take~$L>R$ and require the restriction of a variation to~$B_L$ to be orthogonal to the finitely many negative Neumann eigenfunctions of~$-\Delta-g'(u)$ there.
The interior form is then nonnegative and the exterior contributes only Dirichlet energy.
Finally~$D^2u(0)=-g(a)I/2$ has rank two, so these solutions depend on both coordinates.
\end{proof}

\begin{rmk}
Suppose~$g(s)=f(s_+)$ belongs to~$C^1(\R)$ and that, for some~$p>1$ and positive constants~$c_1,c_2,\delta$,
\begin{align}
 c_1t^p\le f(t)\le c_2t^p,\qquad tf'(t)\ge(1+\delta)f(t)\qquad(t>0).
\end{align}
Put~$w=u_+$ and choose~$\eta$ supported in the region of stability.
Testing the equation with~$w\eta^2$ and the stability inequality with~$w\eta$, as in Lemma~\ref{lemma:annuli}, yields
\begin{align}
 \int_{\R^2}|\nabla(w\eta)|^2\dd x
 &=\int_{\R^2}wf(w)\eta^2\dd x+\int_{\R^2}w^2|\nabla\eta|^2\dd x,\\
 \int_{\R^2} w^2 f'(w)\eta^2\dd x&\le\int_{\R^2}|\nabla(w\eta)|^2\dd x.
\end{align}
Subtracting the common term and using the two hypotheses on~$f$, we obtain
\begin{align}
 \delta c_1\int_{\R^2} w^{p+1}\eta^2\dd x
 &\le\delta\int_{\R^2}wf(w)\eta^2\dd x\\
 &\le\int_{\R^2} w\bigl(wf'(w)-f(w)\bigr)\eta^2\dd x
 \le\int_{\R^2} w^2|\nabla\eta|^2\dd x.
\end{align}
This is the estimate used in Lemma~\ref{lemma:annuli}, with~$\delta c_1$ in place of~$p-1$.
Using the same cutoff and summing over the exterior annuli yields
\begin{align}
 \int_{\R^2\setminus B_R} f(w)\dd x
 \le c_2\int_{\R^2\setminus B_R}w^p\dd x
 \le C R^{-2/(p-1)} .
\end{align}
The integral over~$B_R$ is finite by continuity.
\end{rmk}

\section{Semi-periodic solutions}\label{sec:unified}

We prove Theorem~\ref{thm:main} on~$\R^k\times\R^l$, with~$x\in\R^k$ the transverse variable and~$y\in\R^l$ the periodic variable.
Write~$n=k+l$, where~$k,l\ge1$, and assume
\begin{align}\label{eq:range}
 p>1,\qquad p<\frac{k+2}{k-2}\quad\text{if }k\ge3.
\end{align}
The upper bound involves only the dimension of the radial profile.
In this section we use the same notation for the~$p$-mass on~$\R^n$,
\begin{align}
 M_p(u)=\int_{\R^k\times\R^l}u_+^p\dd x\dd y .
\end{align}

Let~$U$ be the radial solution of~$-\Delta_x U=U_+^p$ on~$\R^k$ normalized by~$U(0)=1$.
Writing~$r=|x|$, its radial profile satisfies
\begin{align}
 -U''(r)-\frac{k-1}{r}U'(r)=U_+(r)^p \qquad (r>0),
\end{align}
\begin{align}
 U(0)=1,\qquad U'(0)=0.
\end{align}
For~$k=1$ we recover the even strip profile~$W_1$.
We write~$B_R=B_R(0)\subset\R^k$ for a ball of radius~$R>0$ and denote the first zero of~$U$ by~$R_*$, so the plasma region of~$U$ is~$B_{R_*}$.
We will describe~$U$ in terms of the positive Dirichlet solution on the unit ball in Section~\ref{subsec:transverse}, where we also study the operator
\begin{align}
 A=-\Delta_x-V,\qquad V=pU_+^{p-1},
\end{align}
on~$L^2_{\mathrm{rad}}(\R^k)$.
By Proposition~\ref{prop:profile}, its lowest eigenvalue~$\lambda_0$ is negative and simple.
Let~$\psi_0>0$ be the corresponding eigenfunction of unit~$L^2$ norm and set~$\omega_0=\sqrt{-\lambda_0}$.

With~$\T=\R/(2\pi\mathbb{Z})$ and~$t=\omega y$, we work on~$\R^k\times\T^l$.
We call a function on the cylinder symmetric if it is radial in~$x$, even in each~$t_j$ and invariant under permutations of the periodic coordinates.
Among the first Fourier modes these symmetries select the one-dimensional space spanned by~$\sum_{j=1}^l\cos t_j$.
Let~$P_0$ denote averaging over the periodic variables,
\begin{align}
 (P_0 v)(x)=\frac{1}{(2\pi)^l}\int_{[0,2\pi]^l}v(x,t)\,\dd t,
\end{align}
where~$t=(t_1,\ldots,t_l)$ and~$\dd t=\dd t_1\cdots\dd t_l$ is product Lebesgue measure on the period cell.
Write~$v=v_0+w$, where~$v_0=P_0v$ and~$P_0w=0$.
To see why we treat the two components separately, fix~$R>R_*$ and suppose that~$h$ is radial in~$x$, has at most polynomial growth and satisfies
\begin{align}
 -\Delta_x h-\omega_0^2\Delta_t h=0\qquad\text{in }(\R^k\setminus\overline{B_R})\times\T^l.
\end{align}
For the bifurcating solutions we will take~$R$ large enough that~$U$ and~$u_s$ are negative for~$|x|\ge R$ and all sufficiently small~$|s|$, as in~\eqref{eq:cylinder}.
For~$|x|>R$ the Fourier expansion in the periodic variables is
\begin{align}\label{eq:exterior-fourier}
 h(x,t)&=h_0(x)+\sum_{\nu\in\mathbb{Z}^l\setminus\{0\}}h_\nu(x)e^{i\nu\cdot t},\\
 h_0(x)&=(P_0h)(x),\qquad
 h_\nu(x)=\frac{1}{(2\pi)^l}\int_{\T^l}h(x,t)e^{-i\nu\cdot t}\dd t.
\end{align}
Since~$-\Delta_t e^{i\nu\cdot t}=|\nu|^2e^{i\nu\cdot t}$, the coefficients satisfy
\begin{align}
 -\Delta_x h_0&=0,\\
 (-\Delta_x+\omega_0^2|\nu|^2)h_\nu&=0,\qquad \nu\in\mathbb{Z}^l\setminus\{0\}.
\end{align}
The radial mean satisfies~$(r^{k-1}h_0')'=0$, so it is affine in~$r$ for~$k=1$, affine in~$\log r$ for~$k=2$ and affine in~$r^{2-k}$ for~$k\ge3$.
For the nonzero modes fix~$0<\beta<\omega_0$ and increase~$R$ so that, for~$r\ge R$ and~$|\nu|\ge1$,
\begin{align}
 (-\Delta_x+\omega_0^2|\nu|^2)e^{\pm\beta r}
 =\left(\omega_0^2|\nu|^2-\beta^2\mp\frac{(k-1)\beta}{r}\right)e^{\pm\beta r}\ge0 .
\end{align}
On an annulus with inner radius~$R$, compare the real and imaginary parts of~$\pm h_\nu$ with
\begin{align}
 |h_\nu(R)|e^{-\beta(r-R)}+\epsilon e^{\beta(r-R)} .
\end{align}
The second term dominates the polynomial bound at every sufficiently large outer radius, so the maximum principle applies.
Letting the outer radius tend to infinity and then~$\epsilon$ to zero proves exponential decay of each mode.
Since~$h(R,\cdot)$ is smooth on~$\T^l$, its Fourier coefficients are absolutely summable, and the same decay holds for~$h-h_0$.
For the explicit expansion on a flat cylinder see~\cite[Section~3, (3.1)--(3.2)]{BMMS2021}, and for the corresponding decomposition on a product with a closed manifold see~\cite[Section~2.2, proofs of Lemmas~2.2 and~2.3]{HNS2024}.
For the inhomogeneous equation used here, Lemma~\ref{lemma:mean} identifies the harmonic part of the mean up to an exponentially decaying remainder, and Lemma~\ref{lemma:free} proves the weighted estimate for the component of zero mean.

Fix~$q>n$ and choose a decay rate~$0<\sigma<\omega_0$.
Let~$\rho$ be a smooth positive radial function on~$\R^k$ such that~$\rho(x)=|x|$ for~$|x|\ge1$, and put~$m_\sigma(x)=e^{\sigma\rho(x)}$.
For~$\Omega=\R^k$ or~$\Omega=\R^k\times\T^l$ define
\begin{align}
 L^q_\sigma(\Omega)
 &=\{h:m_\sigma h\in L^q(\Omega)\},
 &\|h\|_{L^q_\sigma(\Omega)}
 &=\|m_\sigma h\|_{L^q(\Omega)},\\
 W^{2,q}_\sigma(\Omega)
 &=\{v:m_\sigma v\in W^{2,q}(\Omega)\},
 &\|v\|_{W^{2,q}_\sigma(\Omega)}
 &=\|m_\sigma v\|_{W^{2,q}(\Omega)}.
\end{align}
On the cylinder we use the product measure~$\dd x\,\dd t$ and take weak derivatives in both~$x$ and~$t$.
Let~$\mathscr{Y}_0$ be the subspace of radial functions in~$L^q_\sigma(\R^k)$, and set
\begin{align}
 \mathscr{X}_0
 &=\{v_0\in W^{2,q}_{\mathrm{loc,rad}}(\R^k): v_0(0)=0,\ Av_0\in\mathscr{Y}_0\},\\
 \|v_0\|_{\mathscr{X}_0}&=\|Av_0\|_{\mathscr{Y}_0}.
\end{align}
For the component of zero mean we set
\begin{align}
 \mathscr{X}_\perp
 &=\{w\in W^{2,q}_\sigma(\R^k\times\T^l): w\text{ is symmetric and }P_0w=0\},\\
 \mathscr{Y}_\perp
 &=\{h\in L^q_\sigma(\R^k\times\T^l): h\text{ is symmetric and }P_0h=0\}.
\end{align}
We equip~$\mathscr{X}_\perp$ and~$\mathscr{Y}_\perp$ with the inherited norms and the direct sums
\begin{align}\label{eq:spaces}
 \mathscr{X}=\mathscr{X}_0\oplus\mathscr{X}_\perp,\qquad
 \mathscr{Y}=\mathscr{Y}_0\oplus\mathscr{Y}_\perp
\end{align}
with their product norms.

Following the ansatz of Dancer~\cite[Section~1]{Dancer2001} for one periodic variable, we perturb the transverse profile~$U$ by a function periodic in the rescaled variable~$t=\omega y$ and allow the frequency~$\omega$ to vary along the branch.
For~$k=l=1$, the leading perturbation of~$W_1$ is in the direction~$\psi_0(x)\cos t$.
\begin{thm}\label{thm:semi}
Under~\eqref{eq:range}, there exists a continuous branch~$(\omega(s),v_s)$ in~$(0,\infty)\times\mathscr{X}$, parametrized by~$s\in(-\epsilon,\epsilon)$ for some~$\epsilon>0$ and passing through~$(\omega_0,0)$ at~$s=0$, such that the function
\begin{align}
 u_s(x,y)=U(|x|)+v_s(x,\omega(s)y)
\end{align}
satisfies
\begin{align}
 -\Delta u_s=(u_s)_+^p\qquad\text{in }\R^n.
\end{align}
As~$s\to0$ we have
\begin{align}\label{eq:branch}
 v_s(x,t)=s\psi_0(|x|)\sum_{j=1}^l\cos t_j+O_{\mathscr{X}}\bigl(|s|^{\min\{p,2\}}\bigr).
\end{align}
The solutions are radial in~$x$, even in each~$y_j$, invariant under permutations of the~$y_j$ and periodic in each~$y_j$ with common period~$2\pi/\omega(s)$.
For all sufficiently small~$s\ne0$ the positive part is uniformly bounded, its support lies in a fixed transverse cylinder, and~$u_s$ has no nonzero direction of translation invariance.
Each such solution has infinite~$p$-mass and infinite Morse index on~$\R^n$.
\end{thm}

\subsection{The transverse profile}\label{subsec:transverse}
Under~\eqref{eq:range}, let~$\phi_k$ be the unique positive solution of the Dirichlet problem
\begin{align}
\begin{cases}
    -\Delta_x\phi_k=\phi_k^p &\text{in }B_1\subset\R^k,\\
 \phi_k=0 &\text{on }\partial B_1.
\end{cases}
\end{align}
Existence follows by minimizing the Dirichlet energy at fixed~$L^{p+1}$ norm, and every positive solution is radial and strictly decreasing in~$r$, with uniqueness on a fixed ball by scaling~\cite[Theorem~1 and Lemma~2.3]{GNN1979}.
Rescaling~$\phi_k$ to have value one at the origin and matching the radial harmonic extension at the first zero, we obtain
\begin{align}\label{eq:radialtail}
 U(r)=\frac{1}{\phi_k(0)}
 \begin{cases}
  \phi_k(r/R_*), & 0\le r\le R_*,\\[2pt]
  \phi_k'(1)\log(r/R_*), & r>R_*,\ k=2,\\[2pt]
  \displaystyle \phi_k'(1)\frac{(r/R_*)^{2-k}-1}{2-k}, & r>R_*,\ k\ne2 .
 \end{cases}
 \qquad R_*=\phi_k(0)^{(p-1)/2}.
\end{align}
Indeed, the choice of~$R_*$ preserves the radial equation on~$B_{R_*}$, and the exterior formula follows by integrating~$(r^{k-1}U')'=0$ with the same value and first derivative at~$R_*$.
Since~$\phi_k'(1)<0$, the extension is negative for~$r>R_*$, and the radial equation ensures that the second derivatives also agree at~$R_*$.
For~$k\ge3$ the upper bound in~\eqref{eq:range} is necessary for a finite first zero, since the left-hand side of the Pohozaev identity
\begin{align}
 \frac{R_*^k}{2}U'(R_*)^2
 =\left(\frac{k}{p+1}-\frac{k-2}{2}\right)\int_0^{R_*} r^{k-1}U(r)^{p+1}\dd r
\end{align}
is strictly positive.

For~$k=1$, formula~\eqref{eq:quadrature} determines~$W_1$ and its first zero.
For~$k=2$, we have~$\phi_2=\phi$, and~\eqref{eq:radialtail} is the centered profile~$U_{R_*,0}$ in~\eqref{eq:diskprofile}.
For~$k\ge3$ and~$1<p<k/(k-2)$, let~$w_0$ be the ground state centered at the origin in~\cite[Section~6]{BJW2025}, and set~$m=w_0(0)-1>0$.
Then
\begin{align}
 -\Delta_x w_0&=(w_0-1)_+^p,
 &U(r)&=\frac{w_0(m^{-(p-1)/2}r)-1}{m},
\end{align}
since the rescaled function solves the same radial equation and has value one and derivative zero at the origin.

\begin{prop}\label{prop:profile}
Under~\eqref{eq:range}, the operator~$A=-\Delta_x-V$ on~$L^2_{\mathrm{rad}}(\R^k)$, with~$V=pU_+^{p-1}$, has a simple lowest eigenvalue~$\lambda_0<0$.
Its eigenfunction~$\psi_0$ may be chosen positive with unit~$L^2$ norm.
\end{prop}
\begin{proof}
Since~$V$ is bounded and supported in~$\overline{B_{R_*}}$, the compact embedding of~$H^2(B_{R_*})$ into~$L^2(B_{R_*})$ makes multiplication by~$V$ relatively compact with respect to~$-\Delta_x$.
Weyl's theorem on the invariance of the essential spectrum under relatively compact perturbations~\cite[Theorem~6.19]{Teschl2014} thus implies
\begin{align}
 \sigma_{\mathrm{ess}}(A)=\sigma_{\mathrm{ess}}(-\Delta_x)=[0,\infty).
\end{align}
Testing the equation for~$U$ with~$U_+$, we find
\begin{align}
 \int_{\R^k}\bigl(|\nabla U_+|^2-VU_+^2\bigr)\dd x
 &=(1-p)\int_{B_{R_*}} U^{p+1}\dd x<0 .
\end{align}
The bottom of the spectrum is thus a negative eigenvalue.
Taking the absolute value preserves the Rayleigh quotient, and the strong maximum principle makes the resulting normalized eigenfunction~$\psi_0$ strictly positive.
Every regular radial eigenfunction has derivative zero at the origin and is determined by its value there by uniqueness for the radial initial value problem, so the lowest eigenvalue is simple.
\end{proof}

For the weighted space used in Proposition~\ref{prop:linear}, we also need the decay of~$\psi_0$ at infinity.
Outside~$B_{R_*}$ the eigenvalue equation becomes
\begin{align}
 (-\Delta_x+\omega_0^2)\psi_0=0 .
\end{align}
Interior estimates and square integrability imply that~$\psi_0$ tends to zero at infinity.
For~$\sigma<\beta<\omega_0$, the function~$e^{-\beta r}$ is a supersolution of the exterior equation, so comparison followed by interior estimates shows that~$\psi_0$ and its derivatives up to order two are~$O(e^{-\beta r})$.

\subsection{The linearized equation on the cylinder}
\begin{lemma}\label{lemma:mean}
The map~$A:\mathscr{X}_0\to\mathscr{Y}_0$ is an isometric isomorphism.
For every~$R>0$, there exists~$C_R>0$ such that every function~$v_0\in\mathscr{X}_0$ satisfies
\begin{align}
 \|v_0\|_{W^{2,q}(B_R(0))}\le C_R\|v_0\|_{\mathscr{X}_0},\qquad B_R(0)\subset\R^k.
\end{align}
For every~$v_0\in\mathscr{X}_0$ there are constants~$b,c$ and a function~$e$ on~$[R_*,\infty)$ such that
\begin{align}\label{eq:meantail}
 v_0(r)&=c+b\int_{R_*}^r t^{1-k}\dd t+e(r),\qquad r\ge R_*,
\end{align}
with
\begin{align}
 |c|+|b|+\sup_{r\ge R_*}e^{\sigma r}(|e(r)|+|e'(r)|)&\le C\|v_0\|_{\mathscr{X}_0}.
\end{align}
\end{lemma}
\begin{proof}
For~$h\in\mathscr{Y}_0$, the radial equation~$Av_0=h$ is
\begin{align}
 -(r^{k-1}v_0')'=r^{k-1}(h+Vv_0).
\end{align}
Integrating twice with~$v_0(0)=v_0'(0)=0$, we obtain
\begin{align}
 v_0(r)=-\int_0^r t^{1-k}\int_0^t \tau^{k-1}(h+Vv_0)(\tau)\dd\tau\dd t .
\end{align}
On~$[0,r]$, the contribution of~$Vv_0$ to the double integral is bounded by
\begin{align}
 \|V\|_\infty\|v_0\|_\infty
 \int_0^r t^{1-k}\int_0^t \tau^{k-1}\dd\tau\dd t
 =\frac{r^2}{2k}\|V\|_\infty\|v_0\|_\infty .
\end{align}
For the contribution of~$h$, we apply H\"older's inequality with respect to the radial measure and obtain
\begin{align}
 \int_0^t \tau^{k-1}|h(\tau)|\dd\tau
 &\le C t^{k(1-1/q)}\|h\|_{L^q(B_t)},\\
 \int_0^r t^{1-k}\int_0^t \tau^{k-1}|h(\tau)|\dd\tau\dd t
 &\le C\|h\|_{L^q(B_r)}\int_0^r t^{1-k/q}\dd t
 \le C r^{2-k/q}\|h\|_{L^q(B_r)}.
\end{align}
Choose~$a>0$ so that~$a^2\|V\|_\infty/(2k)\le1/2$.
The integral equation is then a contraction on~$C([0,a])$, and its unique solution satisfies
\begin{align}\label{eq:mean-local-value}
 \|v_0\|_{L^\infty(B_a)}\le C a^{2-k/q}\|h\|_{L^q(B_a)} .
\end{align}
Differentiating the integral equation for~$0<r\le a$, we obtain
\begin{align}
 r^{k-1}v_0'(r)=-\int_0^r t^{k-1}(h+Vv_0)(t)\dd t.
\end{align}
Applying the preceding H\"older estimate to the right-hand side, we find
\begin{align}\label{eq:mean-local-derivative}
 |v_0'(r)|\le C r^{1-k/q}\|h+Vv_0\|_{L^q(B_r)} .
\end{align}
Since~$q>k$, the derivative tends to zero at the origin.

To extend the solution to any radius~$R>a$, put~$Y=(v_0,v_0')^{\mathsf T}$ and use the values at~$r=a$ as initial data for
\begin{align}
 Y'(r)=
 \begin{pmatrix}
  0&1\\
  -V(r)&-(k-1)/r
 \end{pmatrix}Y(r)
 +\begin{pmatrix}0\\-h(r)\end{pmatrix},
 \qquad a\le r\le R .
\end{align}
On~$[a,R]$ the matrix is bounded and~$h\in L^1(a,R)$, since~$r^{k-1}$ is bounded above and below by positive constants on~$[a,R]$.
Since~$a$ is fixed and~$V$ is bounded, applying~\eqref{eq:mean-local-derivative} at~$r=a$ and bounding~$v_0$ on~$B_a$ by~\eqref{eq:mean-local-value}, we obtain
\begin{align}
 |Y(a)|\le |v_0(a)|+|v_0'(a)|\le C_a\|h\|_{L^q(B_a)}.
\end{align}
The linear initial value problem has a unique absolutely continuous solution on~$[a,R]$, and Gr\"onwall's inequality implies
\begin{align}
 \sup_{a\le r\le R}|Y(r)|
 &\le C_{a,R}\left(|Y(a)|+\int_a^R|h(t)|\dd t\right)
 \le C_R\|h\|_{\mathscr{Y}_0} .
\end{align}
Uniqueness makes the solutions for different~$R$ agree on their common intervals.

To verify the equation at the origin, first observe that the flux through~$\partial B_\epsilon$, computed with the outward normal of~$B_\epsilon$, is
\begin{align}
 \int_{\partial B_\epsilon}\partial_\nu v_0\dd S
 &=|\sph^{k-1}|\epsilon^{k-1}v_0'(\epsilon)
 =-\int_{B_\epsilon}(h+Vv_0)\dd x\longrightarrow0 .
\end{align}
For~$k=1$ the same formula sums the two outward derivatives at~$x=\pm\epsilon$ of the even extension.
For~$\zeta\in C_c^\infty(B_R)$, integration by parts on~$B_R\setminus\overline{B_\epsilon}$ gives
\begin{align}
 \int_{B_R\setminus\overline{B_\epsilon}}\nabla v_0\cdot\nabla\zeta\dd x
 &=\int_{B_R\setminus\overline{B_\epsilon}}(h+Vv_0)\zeta\dd x
   -v_0'(\epsilon)\int_{\partial B_\epsilon}\zeta\dd S ,
\end{align}
where the minus sign comes from the inward radial normal on the inner boundary of the annulus.
By the flux identity and local integrability,
\begin{align}
 \left|v_0'(\epsilon)\int_{\partial B_\epsilon}\zeta\dd S\right|
 &\le \|\zeta\|_\infty\,|\sph^{k-1}|\epsilon^{k-1}|v_0'(\epsilon)|
 \le \|\zeta\|_\infty\int_{B_\epsilon}|h+Vv_0|\dd x\longrightarrow0 .
\end{align}
Passing to the limit proves the weak equation on the whole ball,
\begin{align}
 \int_{B_R}\nabla v_0\cdot\nabla\zeta\dd x
 =\int_{B_R}(h+Vv_0)\zeta\dd x .
\end{align}
Interior elliptic estimates on~$B_{R+1}$ and the local bounds above now imply
\begin{align}
 \|v_0\|_{W^{2,q}(B_R)}
 \le C_R\bigl(\|h+Vv_0\|_{L^q(B_{R+1})}+\|v_0\|_{L^q(B_{R+1})}\bigr)
 \le C_R\|h\|_{\mathscr{Y}_0} .
\end{align}
Conversely, a radial function in~$\mathscr{X}_0$ belongs to~$C^1$, with value and radial derivative zero at the origin, hence satisfies the same integral equation with~$h=Av_0$.
Since~$V=0$ on~$[R_*,\infty)$, integration of the radial equation from~$R_*$ to~$r$ gives
\begin{align}
 r^{k-1}v_0'(r)
 &=R_*^{k-1}v_0'(R_*)-\int_{R_*}^r t^{k-1}h(t)\dd t
 =b+\int_r^\infty t^{k-1}h(t)\dd t,\\
 b&=R_*^{k-1}v_0'(R_*)-\int_{R_*}^\infty t^{k-1}h(t)\dd t .
\end{align}
Writing~$q'=q/(q-1)$, we use weighted H\"older's inequality to estimate the integral over~$[r,\infty)$ for large~$r$,
\begin{align}
 \int_r^\infty t^{k-1}|h(t)|\dd t
 \le C\|h\|_{\mathscr{Y}_0}
 \left(\int_r^\infty t^{k-1}e^{-q'\sigma t}\dd t\right)^{1/q'}
 \le C\|h\|_{\mathscr{Y}_0}\,r^{(k-1)/q'}e^{-\sigma r}.
\end{align}
For the last bound we substitute~$t=r+s$ and use~$(r+s)^{k-1}\le r^{k-1}(1+s)^{k-1}$ for~$r\ge1$.
It follows that~$v_0'(r)-br^{1-k}=O(e^{-\sigma r}\|h\|_{\mathscr{Y}_0})$.
Integrating this difference from~$r$ to infinity proves~\eqref{eq:meantail}, with the local estimate at~$R_*$ controlling the constants~$b$ and~$c$.
Surjectivity and uniqueness also prove completeness of~$\mathscr{X}_0$.
\end{proof}

\begin{lemma}\label{lemma:free}
The operator~$B_0=-\Delta_x-\omega_0^2\Delta_t$ maps~$\mathscr{X}_\perp$ isomorphically onto~$\mathscr{Y}_\perp$.
\end{lemma}
\begin{proof}
We construct the Green operator for~$B_0$ on functions of zero periodic mean by integrating the heat flow.
The Euclidean and torus Laplacians act on different variables, so their heat operators commute and
\begin{align}
 S_\tau=e^{\tau\Delta_x}e^{\tau\omega_0^2\Delta_t},
 \qquad \partial_\tau S_\tau g=-B_0S_\tau g.
\end{align}
For smooth~$g$ of compact support in~$x$, we integrate in~$\tau$ to obtain
\begin{align}\label{eq:heat-integral-free}
 B_0\int_0^T S_\tau g\,\dd\tau=g-S_Tg.
\end{align}
To pass to the limit as~$T\to\infty$ in~\eqref{eq:heat-integral-free}, we estimate the torus heat operator on functions of zero mean and the Euclidean heat operator in the weighted space.

After removal of the constant Fourier mode, the torus heat kernel is
\begin{align}
 K_\tau(t)=\frac{1}{(2\pi)^l}
 \sum_{\nu\in\mathbb{Z}^l\setminus\{0\}}
 e^{-\omega_0^2\tau|\nu|^2}e^{i\nu\cdot t}.
\end{align}
For~$\tau\ge1$, each index in the sum satisfies~$|\nu|^2\ge1$, whence
\begin{align}
 \int_{\T^l}|K_\tau(t)|\,\dd t
 \le e^{-\omega_0^2\tau}
 \sum_{\nu\in\mathbb{Z}^l\setminus\{0\}}e^{-\omega_0^2(|\nu|^2-1)}.
\end{align}
The sum is finite.
For~$0<\tau\le1$, we use the full torus heat kernel, which is the periodization of the Euclidean Gaussian,
\begin{align}
 H_\tau^{\T}(t)
 =\frac{1}{(4\pi\omega_0^2\tau)^{l/2}}
   \sum_{m\in\mathbb Z^l}
   \exp\left(-\frac{|t+2\pi m|^2}{4\omega_0^2\tau}\right).
\end{align}
Every summand is real and positive, and integration over one period cell gives
\begin{align}
 \int_{\T^l} H_\tau^{\T}(t)\dd t
 =\frac{1}{(4\pi\omega_0^2\tau)^{l/2}}
   \int_{\R^l}e^{-|t|^2/(4\omega_0^2\tau)}\dd t=1.
\end{align}
The kernel~$K_\tau=H_\tau^{\T}-(2\pi)^{-l}$ changes sign, since its integral is zero and its Fourier series is positive at~$t=0$.
For~$P_0g=0$, the constant term contributes nothing to convolution in~$t$, so Young's inequality implies
\begin{align}
 \|K_\tau*g(x,\cdot)\|_{L^q(\T^l)}
 =\|H_\tau^{\T}*g(x,\cdot)\|_{L^q(\T^l)}
 \le\|g(x,\cdot)\|_{L^q(\T^l)} .
\end{align}
Combining the bounds for small and large times, the torus heat operator on functions of zero mean has norm at most~$Ce^{-\omega_0^2\tau}$ for all~$\tau>0$.

For the Euclidean factor, the weight~$m_\sigma$ is comparable to~$e^{\sigma|x|}$, so the triangle inequality implies
\begin{align}
 m_\sigma(x)\le Ce^{\sigma|z|}m_\sigma(x-z).
\end{align}
With~$H_\tau(z)=(4\pi\tau)^{-k/2}e^{-|z|^2/(4\tau)}$, we pass to polar coordinates and put~$r=2\sqrt{\tau}\,s$, which leads to
\begin{align}
 \int_{\R^k}H_\tau(z)e^{\sigma|z|}\,\dd z
 &=C_k\int_0^\infty s^{k-1}e^{-s^2+2\sigma\sqrt{\tau}\,s}\,\dd s\\
 &=C_ke^{\sigma^2\tau}\int_0^\infty s^{k-1}e^{-(s-\sigma\sqrt{\tau})^2}\,\dd s\\
 &\le C(1+\tau^{(k-1)/2})e^{\sigma^2\tau}.
\end{align}
For the last inequality, put~$a=\sigma\sqrt{\tau}$ and~$b=s-a$, and bound~$(a+b)^{k-1}$ on~$b\ge-a$ by~$C(1+a)^{k-1}(1+|b|)^{k-1}$.
The latter factor is integrable against~$e^{-b^2}\,\dd b$ on~$\R$.

By Young's inequality for the Euclidean convolution and the torus bound, for~$P_0g=0$ we have
\begin{align}
 \|S_\tau g\|_{L^q_\sigma}
 &\le C(1+\tau^{(k-1)/2})e^{-(\omega_0^2-\sigma^2)\tau}\|g\|_{L^q_\sigma}.
\end{align}
Since~$\sigma<\omega_0$, the right-hand side is integrable in~$\tau$, so we may define~$\mathcal{G}g$ by a convergent integral in~$L^q_\sigma$, with
\begin{align}
 \mathcal{G}g=\int_0^\infty S_\tau g\,\dd\tau,
 \qquad \|\mathcal{G}g\|_{L^q_\sigma}\le C\|g\|_{L^q_\sigma}.
\end{align}
Letting~$T\to\infty$ in~\eqref{eq:heat-integral-free} and then approximating~$g$ in~$L^q_\sigma$ proves~$B_0\mathcal{G}g=g$ in distributions.

To estimate the derivatives of~$w=\mathcal{G}g$, cover the cylinder by unit balls~$D$ whose doubled balls~$2D$ have bounded overlap.
The interior elliptic estimate reads
\begin{align}
 \|w\|_{W^{2,q}(D)}
 \le C\bigl(\|g\|_{L^q(2D)}+\|w\|_{L^q(2D)}\bigr).
\end{align}
On~$2D$ the weight is comparable to its value at the center, and its first two derivatives are bounded by~$Cm_\sigma$.
Multiplying by the weight at each center and summing the~$q$-th powers, we obtain
\begin{align}
 \|w\|_{W^{2,q}_\sigma}
 \le C\bigl(\|g\|_{L^q_\sigma}+\|w\|_{L^q_\sigma}\bigr)
 \le C\|g\|_{L^q_\sigma}.
\end{align}
The heat flow preserves the mean and the imposed symmetries, hence~$\mathcal{G}$ maps~$\mathscr{Y}_\perp$ into~$\mathscr{X}_\perp$.
For smooth~$v\in\mathscr{X}_\perp$ of compact support in~$x$,~$S_\tau B_0v=B_0S_\tau v$, so~\eqref{eq:heat-integral-free} also implies~$\mathcal{G}B_0v=v$.
Radial cutoffs in~$x$ and symmetric mollification approximate any~$v\in\mathscr{X}_\perp$ in~$W^{2,q}_\sigma$ while preserving its symmetries and zero mean.
Thus~$\mathcal{G}B_0v=v$ throughout~$\mathscr{X}_\perp$.
\end{proof}

In the planar case~$k=l=1$, the operator~$B_0$ on the~$j$-th cosine mode is~$-\partial_x^2+a^2$, where~$a=j\omega_0$ and~$j\ge1$.
Its Green function on the line satisfies
\begin{align}
 G_a(x)&=\frac{e^{-a|x|}}{2a},\\
 \int_{\R} G_a(x)e^{\sigma|x|}\dd x
 &=\frac1a\int_0^\infty e^{-(a-\sigma)x}\dd x=\frac{1}{a(a-\sigma)}
 \qquad (0<\sigma<a).
\end{align}
The first cosine mode thus determines the common decay threshold~$\omega_0$.
The separation into a transverse eigenfunction and a first periodic mode is the same as in the local bifurcation analysis of~\cite[Section~1]{Dancer2001}.

For each~$\omega>0$, define the bounded linear operator
\begin{align}
 L_\omega\colon\mathscr{X}&\longrightarrow\mathscr{Y},\\
 L_\omega v &=(A-\omega^2\Delta_t)v
 =-\Delta_x v-\omega^2\Delta_t v-pU_+^{p-1}v .
\end{align}
We now identify its kernel at the bifurcation frequency.
On~$\Omega=\R^k\times\T^l$, set
\begin{align}
 \Phi(x,t)=\psi_0(|x|)\sum_{j=1}^l\cos t_j .
\end{align}
Define the continuous linear functional~$\ell\colon\mathscr{Y}\to\R$ by
\begin{align}
 \ell(g)=\int_{\R^k\times\T^l}g(x,t)\Phi(x,t)\dd x\dd t .
\end{align}

\begin{prop}\label{prop:linear}
At~$\omega=\omega_0$ the kernel and range are
\begin{align}\label{eq:kernel}
 \ker L_{\omega_0}=\R\Phi,\qquad
 \operatorname{Ran} L_{\omega_0}=\ker\ell .
\end{align}
\end{prop}
\begin{proof}
Since each function in~$\mathscr{X}_0$ is independent of~$t$, the restriction of~$L_{\omega_0}$ to~$\mathscr{X}_0$ coincides with~$A$, which maps~$\mathscr{X}_0$ isomorphically onto~$\mathscr{Y}_0$ by Lemma~\ref{lemma:mean}.
Multiplication by~$V$ maps~$\mathscr{X}_\perp$ compactly into~$\mathscr{Y}_\perp$ by its compact support in~$x$ and local Rellich compactness.
By Lemma~\ref{lemma:free}, the restriction of~$L_{\omega_0}=B_0-V$ to~$\mathscr{X}_\perp$ is a compact perturbation of an isomorphism onto~$\mathscr{Y}_\perp$.
Since~$L_{\omega_0}$ preserves the mean and zero-mean subspaces, it is Fredholm of index zero on their direct sum.
If~$L_{\omega_0}w=0$, its mean satisfies~$A(P_0w)=0$, hence~$P_0w=0$.
For each~$\nu\in\mathbb{Z}^l\setminus\{0\}$, write the Fourier coefficient of~$w$ as in~\eqref{eq:exterior-fourier},
\begin{align}
 w_\nu(x)=\frac{1}{(2\pi)^l}\int_{\T^l}w(x,t)e^{-i\nu\cdot t}\dd t .
\end{align}
The operator~$A$ acts only on~$x$, so the equation~$Aw=\omega_0^2\Delta_t w$ and integration by parts on the torus imply
\begin{align}
 Aw_\nu(x)
 &=\frac{1}{(2\pi)^l}\int_{\T^l}(Aw)(x,t)e^{-i\nu\cdot t}\dd t\\
 &=\frac{\omega_0^2}{(2\pi)^l}\int_{\T^l}(\Delta_t w)(x,t)e^{-i\nu\cdot t}\dd t\\
 &=\frac{\omega_0^2}{(2\pi)^l}\int_{\T^l}w(x,t)\Delta_t(e^{-i\nu\cdot t})\dd t
 =-\omega_0^2|\nu|^2w_\nu(x).
\end{align}
Each coefficient belongs to~$H^2_{\mathrm{rad}}(\R^k)$ by the exponential weight and~$q>2$.
Testing against~$\overline{w_\nu}$, recalling~$\lambda_0=-\omega_0^2$ and using the variational characterization of the lowest radial eigenvalue, we find
\begin{align}
 \lambda_0|\nu|^2\|w_\nu\|_2^2
 =\int_{\R^k}\bigl(|\nabla w_\nu|^2-V|w_\nu|^2\bigr)\dd x
 \ge\lambda_0\|w_\nu\|_2^2.
\end{align}
Since~$\lambda_0<0$, this forces~$w_\nu=0$ for~$|\nu|^2>1$.
At~$|\nu|=1$ each radial coefficient is a multiple of~$\psi_0$, evenness removes the sine modes, and the permutation symmetry makes the cosine coefficients equal.
The exponential decay of~$\psi_0$ and its first two derivatives ensures~$\Phi\in\mathscr{X}_\perp$.
Moreover,
\begin{align}
 A\Phi=\lambda_0\Phi,\qquad \Delta_t\Phi=-\Phi,\qquad
 L_{\omega_0}\Phi=(\lambda_0+\omega_0^2)\Phi=0 .
\end{align}
Hence the kernel is spanned by~$\Phi$.

For the range, write~$v=v_0+w$ with~$v_0=P_0v$ and~$w\in\mathscr{X}_\perp$.
Since~$\Phi$ has zero mean in~$t$, we have
\begin{align}
 \ell(L_{\omega_0}v_0)
 =\int_{\R^k}(Av_0)(x)\left(\int_{\T^l}\Phi(x,t)\dd t\right)\dd x=0 .
\end{align}
The functions~$w$ and~$\Phi$ are periodic in~$t$ and decay exponentially in~$x$ together with their first derivatives, so integration by parts gives
\begin{align}
 \ell(L_{\omega_0}v)
 &=\int_{\Omega}(L_{\omega_0}w)\Phi\dd x\dd t\\
 &=\int_{\Omega}\bigl(\nabla_xw\cdot\nabla_x\Phi+
   \omega_0^2\nabla_tw\cdot\nabla_t\Phi-Vw\Phi\bigr)\dd x\dd t\\
 &=\int_{\Omega}w\bigl(A\Phi-\omega_0^2\Delta_t\Phi\bigr)\dd x\dd t
 =(\lambda_0+\omega_0^2)\int_{\Omega}w\Phi\dd x\dd t=0 .
\end{align}
The range is thus contained in~$\ker\ell$.
Its codimension is one by the Fredholm index and the kernel computation, and~$\ker\ell$ also has codimension one since~$\ell(\Phi)>0$.
This proves~\eqref{eq:kernel}.
\end{proof}

\subsection{Bifurcation from the radial solution}\label{subsec:bifurcation}

Let~$g(s)=s_+^p$ for~$s\in\R$.
At~$U$ the linearized potential~$V$ and the nonlinear remainder~$N$ are
\begin{align}
 V&=g'(U)=pU_+^{p-1},\\
 N(v)&=g(U+v)-g(U)-Vv.
\end{align}
We seek a solution of the plasma equation in the form
\begin{align}
 -\Delta_x u-\Delta_y u&=g(u)\qquad\text{in }\R^k\times\R^l,\\
 u(x,y)&=U(|x|)+v(x,\omega y).
\end{align}
Setting~$t=\omega y$ and applying the chain rule in~$y$, we obtain
\begin{align}
 -\Delta_x U-\Delta_x v-\omega^2\Delta_t v=g(U+v).
\end{align}
Subtracting the profile equation~$-\Delta_xU=g(U)$ and then subtracting~$Vv$ from both sides, we obtain the equivalent equation on~$\R^k\times\T^l$,
\begin{align}\label{eq:perturbation}
 L_\omega v-N(v)=0.
\end{align}
To parametrize the branch by the coefficient of~$\Phi$, define the continuous linear functional
\begin{align}
 J\colon\mathscr{X}\longrightarrow\R,\qquad
 J(v)=\frac{\ell((I-P_0)v)}{\ell(\Phi)} .
\end{align}
Since~$\psi_0$ has unit~$L^2(\R^k)$ norm and the functions~$\cos t_j$ are pairwise orthogonal on~$\T^l$, we have
\begin{align}
 \ell(\Phi)
 &=\int_{\R^k}\psi_0(|x|)^2\dd x\sum_{j=1}^l\int_{\T^l}\cos^2 t_j\dd t
 =\frac{l}{2}(2\pi)^l .
\end{align}
In particular,~$\ell(\Phi)=\pi$ in the planar case, and~$J(\Phi)=1$ since~$P_0\Phi=0$.
By Proposition~\ref{prop:linear}, the restriction of~$L_{\omega_0}$ to~$\mathscr{X}_1=\ker J$ is an isomorphism onto~$\ker\ell$.

\begin{prop}\label{prop:branch}
For all sufficiently small~$|s|$ there is a continuous branch~$(\omega(s),v_s)$ of solutions of~\eqref{eq:perturbation} with~$J(v_s)=s$, satisfying~\eqref{eq:branch} and
\begin{align}
 \omega(s)=\omega_0+O\bigl(|s|^{\min\{p-1,1\}}\bigr).
\end{align}
If~$p>2$, the branch is of class~$C^1$.
\end{prop}
\begin{proof}
First we choose~$R>R_*$ so that~$U+v<0$ for~$|x|\ge R$ whenever~$\|v\|_{\mathscr{X}}$ is sufficiently small.
Indeed, write~$v=v_0+w$ with~$P_0w=0$.
For~$k=1,2$, the coefficient of the affine or logarithmic term of~$U$ in~\eqref{eq:radialtail} is strictly negative, and by~\eqref{eq:meantail} the corresponding coefficient of~$v_0$ tends to zero with~$\|v\|_{\mathscr{X}}$.
For~$k\ge3$, the limit of~$U$ at infinity is strictly negative and~\eqref{eq:meantail} bounds~$v_0$ uniformly on~$[R_*,\infty)$.
In every dimension Sobolev embedding implies
\begin{align}
 |w(x,t)|\le C e^{-\sigma|x|}\|w\|_{W^{2,q}_\sigma}.
\end{align}
These bounds allow the same choice of~$R$ for all sufficiently small perturbations.

On~$\overline{B_R}\times\T^l$, the weight~$m_\sigma$ is bounded above and below by positive constants and Lemma~\ref{lemma:mean} together with Sobolev embedding implies
\begin{align}
 \|v\|_{C^0(\overline{B_R}\times\T^l)}\le C_R\|v\|_{\mathscr{X}} .
\end{align}
For sufficiently small~$\|v\|_{\mathscr{X}}$, the values of~$U$ and~$U+v$ on~$\overline{B_R}\times\T^l$ lie in a fixed bounded interval on which~$g'$ is H\"older continuous with exponent~$\alpha=\min\{p-1,1\}$.
The integral formula for the remainder and the expression for its derivative imply
\begin{align}
 |N(v)|&=\left|\int_0^1\bigl(g'(U+\theta v)-g'(U)\bigr)v\dd\theta\right|
 \le C|v|^{1+\alpha},\\
 |DN(v)h|&=|\bigl(g'(U+v)-g'(U)\bigr)h|
 \le C|v|^\alpha|h|.
\end{align}
Applying the first estimate with base point~$U+v$ bounds the remainder after linearization by~$C|h|^{1+\alpha}$.
For~$|x|\ge R$ both~$N(v)$ and~$DN(v)h$ vanish.
Taking weighted~$L^q$ norms on~$B_R\times\T^l$ and using the H\"older continuity of~$g'$ proves that~$N$ is~$C^1$ near zero, with
\begin{align}\label{eq:remainder}
 \|N(v)\|_{\mathscr{Y}}\le C\|v\|_{\mathscr{X}}^{1+\alpha},\qquad
 \|DN(v)\|_{\mathscr{X}\to\mathscr{Y}}\le C\|v\|_{\mathscr{X}}^\alpha .
\end{align}
When~$p>2$, continuity of~$g''$ on the bounded range of~$U+v$ in~$\overline{B_R}\times\T^l$ likewise implies~$N\in C^2$ near zero.

Now write~$v=s(\Phi+z)$ with~$z\in\mathscr{X}_1$, so that~$J(v)=s$.
Let~$I\subset(0,\infty)$ be an open interval containing~$\omega_0$ and let~$B_\delta^{\mathscr{X}_1}(0)$ denote the ball of radius~$\delta$ in~$\mathscr{X}_1$.
For sufficiently small~$\epsilon,\delta>0$, define
\begin{align}
 E\colon(-\epsilon,\epsilon)\times I\times B_\delta^{\mathscr{X}_1}(0)
 &\longrightarrow\mathscr{Y},\\
 E(s,\omega,z)&=
 \begin{cases}
 L_\omega(\Phi+z)-s^{-1}N(s(\Phi+z)),&s\ne0,\\
 L_\omega(\Phi+z),&s=0 .
 \end{cases}
\end{align}
For~$s\ne0$ we have
\begin{align}
 L_\omega v-N(v)=sE(s,\omega,z),
\end{align}
so each point of~$E^{-1}(0)$ with~$s\ne0$ determines a nontrivial solution of~\eqref{eq:perturbation}.
Conversely, let~$(\omega,v)$ be a solution near~$(\omega_0,0)$, and put~$s=J(v)$ and~$\zeta=v-s\Phi\in\mathscr{X}_1$.
Since~$L_{\omega_0}\Phi=0$, we have
\begin{align}
 L_{\omega_0}\zeta=N(v)+(\omega^2-\omega_0^2)\Delta_t v .
\end{align}
The bounded inverse of~$L_{\omega_0}$ on~$\mathscr{X}_1$ and~\eqref{eq:remainder} imply
\begin{align}
 \|\zeta\|_{\mathscr{X}}
 \le C\bigl(\|v\|_{\mathscr{X}}^\alpha+|\omega^2-\omega_0^2|\bigr)\|v\|_{\mathscr{X}} .
\end{align}
Choosing the neighborhood so that the coefficient of~$\|v\|_{\mathscr{X}}$ is at most~$1/2$, we obtain
\begin{align}
 \|v\|_{\mathscr{X}}\le |s|\|\Phi\|_{\mathscr{X}}+\tfrac12\|v\|_{\mathscr{X}},
 \qquad \|v\|_{\mathscr{X}}\le 2|s|\|\Phi\|_{\mathscr{X}} .
\end{align}
Hence every nearby solution with~$v\ne0$ has~$s\ne0$, and the estimate for~$\zeta$ implies that~$z=\zeta/s$ tends to zero as~$(\omega,v)\to(\omega_0,0)$ along such solutions.
Thus~$E^{-1}(0)$ with~$s\ne0$ describes all nearby nontrivial solutions.
By~\eqref{eq:remainder}, the map~$E$ is continuous and its partial derivatives in~$(\omega,z)$ depend continuously on all three variables.
Indeed, for~$s\ne0$ the derivative in~$z$ of the divided nonlinear term is the restriction of~$DN(s(\Phi+z))$ to~$\mathscr{X}_1$, which tends to zero as~$s\to0$.

Next we show that~$D_{(\omega,z)}E$ at the bifurcation point is invertible, so that we can apply the implicit function theorem with~$s$ as a continuous parameter.
Since~$L_{\omega_0}\Phi=0$ and~$\Delta_t\Phi=-\Phi$, we have~$E(0,\omega_0,0)=0$ and
\begin{align}
 D_{(\omega,z)}E(0,\omega_0,0)\colon\R\times\mathscr{X}_1
 &\longrightarrow\mathscr{Y},\\
 (b,h)&\longmapsto 2\omega_0 b\Phi+L_{\omega_0}h .
\end{align}
Given~$H\in\mathscr{Y}$, applying~$\ell$ to the equation for~$(b,h)$ and using~\eqref{eq:kernel} determines~$b$, after which~$h$ satisfies
\begin{align}
 b&=\frac{\ell(H)}{2\omega_0\ell(\Phi)},\\
 L_{\omega_0}h&=H-\frac{\ell(H)}{\ell(\Phi)}\Phi .
\end{align}
The right-hand side of the second equation lies in~$\ker\ell$, so the isomorphism from~$\mathscr{X}_1$ onto~$\ker\ell$ determines~$h$ uniquely and continuously in~$H$.
The implicit function theorem now supplies unique continuous functions~$\omega(s)$ and~$z(s)$ with~$\omega(0)=\omega_0$ and~$z(0)=0$, and after shrinking the neighborhoods we have
\begin{align}
 E^{-1}(0)=\bigl\{(s,\omega(s),z(s)):|s|<\epsilon\bigr\}.
\end{align}
Returning to~$v_s=s(\Phi+z(s))$, we obtain the bifurcating branch, which meets the trivial branch~$v=0$ at~$(\omega_0,0)$.

For the estimate, put~$X=(\omega,z)$,~$X_*=(\omega_0,0)$ and~$D=D_XE(0,X_*)$.
By continuity of~$D_XE$, the map~$T_s(X)=X-D^{-1}E(s,X)$ satisfies~$\|D_XT_s\|\le1/2$ on a fixed ball about~$X_*$ for all sufficiently small~$|s|$.
Writing~$X_s=(\omega(s),z(s))$ and using~$X_s=T_s(X_s)$, we obtain
\begin{align}
 \|X_s-X_*\|&\le\tfrac12\|X_s-X_*\|+\|D^{-1}E(s,X_*)\|,\\
 |\omega(s)-\omega_0|+\|z(s)\|_{\mathscr{X}}
 &\le C\|E(s,\omega_0,0)\|_{\mathscr{Y}}
 \le C|s|^\alpha
\end{align}
by~\eqref{eq:remainder}.
Multiplication by~$s$ proves~\eqref{eq:branch}.
For~$p>2$, put~$h=\Phi+z$ and write
\begin{align}
 \frac{N(sh)}s=\int_0^1 DN(\theta sh)[h]\dd\theta\qquad(s\ne0).
\end{align}
The integral vanishes at~$s=0$ by~$DN(0)=0$.
Since~$N$ is~$C^2$, the integral is~$C^1$ in~$(s,h)$ and the differentiable implicit function theorem applies, as in simple-eigenvalue bifurcation~\cite{CR}.
\end{proof}

\subsection{The plasma region and the Morse index}\label{subsec:positivity-index}

\begin{proof}[Proof of Theorem~\ref{thm:semi}]
For the branch in Proposition~\ref{prop:branch}, set
\begin{align}
 u_s(x,y)=U(|x|)+v_s(x,\omega(s)y).
\end{align}
The equation holds distributionally, and local Sobolev embedding followed by interior elliptic regularity shows that~$u_s$ is classical.
The symmetries follow from the definition of~$\mathscr{X}$.
By the uniform exterior negativity established in the proof of Proposition~\ref{prop:branch} and the local uniform convergence to~$U$, there exist~$0<r_*<R_*<R$ and~$C>0$ such that
\begin{align}\label{eq:cylinder}
 u_s\ge\tfrac12\quad (|x|\le r_*),\qquad
 (u_s)_+\le C,\qquad
 u_s<0\quad (|x|\ge R).
\end{align}
In particular~$u_s$ and~$U$ are harmonic for~$|x|>R$, since~$R>R_*$, so the exterior equation for the perturbation is exactly the free equation at frequency~$\omega(s)$.
On a period cell~$\mathcal C$ in the~$y$ variables we have
\begin{align}
 0<\int_{\R^k\times\mathcal C}(u_s)_+^p\dd x\dd y<\infty.
\end{align}
Summing over all period cells shows that~$M_p(u_s)=\infty$.

Fix~$s\ne0$ and suppose~$u_s$ is invariant under every translation along~$(a,b)\in\R^k\times\R^l$.
Translating a point of the positive core and using~\eqref{eq:cylinder} shows that~$a=0$.
By permutation symmetry the~$l$ integrals below are equal, and by the definition of~$J$ their sum is~$s\ell(\Phi)$, so for each~$j$ we have
\begin{align}
 \int_{\R^k\times\T^l}v_s(x,t)\psi_0(|x|)\cos t_j\dd x\dd t=\frac{s}{l}\ell(\Phi)\ne0.
\end{align}
If~$v_s$ is invariant under~$t\mapsto t+\tau\omega(s)b$, then a change of variables on the torus multiplies this integral by~$\cos(\tau\omega(s)b_j)$, since the corresponding sine integral vanishes by evenness.
Thus~$\cos(\tau\omega(s)b_j)=1$ for every~$\tau\in\R$, which forces~$b_j=0$.
A function of fewer than~$n$ linear coordinates would have a nonzero direction of invariance.

It remains to compute the full-space index.
Put~$w=(u_s)_+$ and choose a cutoff~$\chi_L(y)$ equal to one on~$B_L^l$, supported in~$B_{2L}^l$, with~$|\nabla\chi_L|\le C/L$.
Testing the original equation with~$w\chi_L^2$ and expanding the gradient as in Lemma~\ref{lemma:annuli}, we obtain
\begin{align}
 Q_{u_s}(w\chi_L)
 =(1-p)\int_{\R^k\times\R^l}w^{p+1}\chi_L^2\dd x\dd y+\int_{\R^k\times\R^l}w^2|\nabla\chi_L|^2\dd x\dd y.
\end{align}
By~\eqref{eq:cylinder},~$w\ge1/2$ on~$B_{r_*}^k\times\R^l$, while~$w$ is bounded and vanishes for~$|x|\ge R$.
Thus
\begin{align}
 \int_{\R^k\times\R^l}w^{p+1}\chi_L^2\dd x\dd y
 &\ge 2^{-(p+1)}|B_{r_*}^k|\,|B_L^l|\ge cL^l,\\
 \int_{\R^k\times\R^l}w^2|\nabla\chi_L|^2\dd x\dd y
 &\le CL^{-2}|B_R^k|\,|B_{2L}^l\setminus B_L^l|\le C'L^{l-2},\\
 Q_{u_s}(w\chi_L)&\le -(p-1)cL^l+C'L^{l-2}.
\end{align}
For large~$L$ this is a compactly supported negative direction, which can be approximated by a smooth one.
Translating it by sufficiently separated integer multiples of a period produces disjoint negative directions in any prescribed number.
This proves the index assertion in every admissible splitting.
Taking~$k=l=1$ proves Theorem~\ref{thm:main}, since~$U$ and~$W_1$ solve the same one-dimensional initial-value problem with value one and derivative zero at the origin.
\end{proof}

\begin{rmk}
On a single period cell, with periodic variations, the index of~$u_s$ is finite and positive.
Write~$T_s=2\pi/\omega(s)$ and choose~$R$ as in~\eqref{eq:cylinder}.
Set
\begin{align}
 \Omega_s=\R^k\times(\R/T_s\mathbb{Z})^l,\qquad K=B_R^k\times(\R/T_s\mathbb{Z})^l.
\end{align}
Since~$w=(u_s)_+$ vanishes outside~$K$, we can test the equation with~$w$ on the full period cell to obtain
\begin{align}
 Q_{u_s}(w)=(1-p)\int_{\Omega_s}w^{p+1}\dd x\dd y<0.
\end{align}
On~$K$ the operator~$-\Delta-p(u_s)_+^{p-1}$ with Neumann condition at~$|x|=R$ has finitely many negative eigenvalues, say~$m$, counted with multiplicity.
If the restriction of a periodic variation~$\zeta$ to~$K$ is orthogonal in~$L^2(K)$ to the corresponding eigenfunctions, then its quadratic form on~$K$ is nonnegative.
Since the potential vanishes outside~$K$, we obtain
\begin{align}
 Q_{u_s}(\zeta)
 &=\int_K\bigl(|\nabla\zeta|^2-p(u_s)_+^{p-1}\zeta^2\bigr)\dd x\dd y
   +\int_{\Omega_s\setminus K}|\nabla\zeta|^2\dd x\dd y\ge0 .
\end{align}
These orthogonality conditions have codimension at most~$m$, so the index is at most~$m$.

In the planar case~$k=l=1$, the branch passes at~$s=0$ through the solution~$u_0(x,y)=W_1(x)$, whose plasma region is~$(-R_*,R_*)\times\R$.
Since~$u_0$ is independent of~$y$, we may compute its Morse index on~$\R\times(\R/T\mathbb{Z})$ for any period~$T>0$.
The transverse operator~$-\partial_x^2-p(W_1)_+^{p-1}$ on the full line has exactly one negative eigenvalue,~$\lambda_0$.
In fact~$(W_1)_+$ is a negative direction, while the ground-state identity with~$|W_1'|$ on each half-line shows that the form is nonnegative on tests vanishing at zero.
The even positive eigenfunction from Proposition~\ref{prop:profile} realizes the negative eigenvalue.
The~$j$-th Fourier mode in~$y$ shifts the transverse eigenvalue by~$(2\pi j/T)^2$, so the cell index is
\begin{align}
 \#\{j\in\mathbb{Z}:\lambda_0+(2\pi j/T)^2<0\}.
\end{align}
At the bifurcation period~$T_0=2\pi/\omega_0$, we have~$\lambda_0+(2\pi/T_0)^2=0$, so the cell index is one and the linearized operator has the two linearly independent zero eigenfunctions
\begin{align}
 \psi_0(x)\cos(\omega_0 y),\qquad \psi_0(x)\sin(\omega_0 y).
\end{align}
Only the cosine mode belongs to the space of functions even in~$y$ used in the bifurcation argument.
\end{rmk}

\ 

\textbf{Declaration of generative AI and AI-assisted technologies in the manuscript preparation process.}
The authors used ChatGPT (OpenAI) to obtain feedback on mathematical arguments and to improve the clarity and presentation of the manuscript.
The authors critically reviewed and revised the resulting material, independently verified the mathematical arguments and references, and accept full responsibility for the accuracy and integrity of the final manuscript.

\medskip

\textbf{Conflict of Interest} The authors have no conflicts to disclose. 

\medskip

\textbf{Data Availability Statement} Data sharing is not applicable to this article as no datasets were generated or analyzed during the current study.

\printbibliography[heading=bibintoc,title={References}]

\end{document}